\documentclass[12pt,fleqn]{article}
\usepackage{amsfonts}
\usepackage{amssymb}
\usepackage{amsmath}
\usepackage{natbib}
\usepackage{amstext,amsthm}
\usepackage{mathtools}
\usepackage{xr}
\usepackage[labelsep=none, skip=-18pt]{caption}

\renewcommand{\mathbf}{\boldsymbol}
\renewcommand{\thepage}{}
\renewcommand{\appendix}{\footnotesize\parindent 0cm\setcounter{equation}{0}
\renewcommand{\theequation}{A.\arabic{equation}}
\setcounter{lemma}{0}\renewcommand{\thelemma}{A.\arabic{lemma}}}

\newcommand{\descr}{\mathrm{descr}}

\newcommand{\causal}{{\rm causal}}
\newcommand{\desc}{{\rm descr}}

\newcommand{\ehw}{{\rm ehw}}

\newcommand{\sample}{{\rm sample}}

\newcommand{\cond}{{\rm cond}}

\newtheorem{assumption}{Assumption}
\newtheorem{theorem}{Theorem}

\newtheorem{lemma}{Lemma}

\newtheorem{definition}{Definition}
\theoremstyle{definition}
\newtheorem{remark}{Comment}
\def\monthname{\ifcase\month\or
  January\or February\or March\or April\or May\or June\or July\or
  August\or September\or October\or November\or December\fi}

\renewcommand{\appendix}{\small\parindent 0cm\setcounter{equation}{0}
\parskip 5pt
\renewcommand{\theequation}{A.\arabic{equation}}
\setcounter{lemma}{0}\renewcommand{\thelemma}{A.\arabic{lemma}}
\setcounter{theorem}{0}\renewcommand{\thetheorem}{A.\arabic{theorem}}}
\numberwithin{equation}{section}
\allowdisplaybreaks

\begin{document}

\title{Sampling-based \textit{vs.}~Design-based Uncertainty in Regression
Analysis \thanks{{\small We are grateful for comments by Daron Acemoglu,
Joshua Angrist, Matias Cattaneo, Alex Lenk, Jonas Metzger, Evan Munro, Michael Pollmann, Jim Poterba, Karthik Rajkumar, Tymon S\l %
{}oczy\'{n}ski, Bas Werker, and seminar participants at Microsoft Research,
Michigan, Brown University, MIT, Stanford, Princeton, NYU, Columbia, Tilburg
University, the Tinbergen Institute, American University, Montreal, Michigan
State, Maryland, Pompeu Fabra, Carlos III, and University College London,
three referees, and especially for discussions with Gary Chamberlain. An
earlier version of this paper circulated under the title ``Finite Population
Causal Standard Errors'' (\citet{abadieathey}). }}}
\author{ Alberto Abadie\thanks{{\small Professor of Economics, Massachusetts
Institute of Technology, and NBER, abadie@mit.edu.}} \and Susan Athey\thanks{%
{\small Professor of Economics, Graduate School of Business, Stanford
University, and NBER, athey@stanford.edu. }} \and Guido W. Imbens\thanks{%
{\small Professor of Economics, Graduate School of Business, and Department
of Economics, Stanford University, and NBER, imbens@stanford.edu. }} \and %
Jeffrey M. Wooldridge\thanks{{\small University Distinguished Professor,
Department of Economics, Michigan State University, wooldri1@msu.edu}} }
\date{ Current version \ifcase\month\or
January\or February\or March\or April\or May\or June\or
July\or August\or September\or October\or November\or December\fi \ \number%
\year\ \ -- First version September 2013 }

\begin{titlepage}
\clearpage\thispagestyle{empty}
\maketitle

\begin{abstract}
\noindent Consider a researcher estimating the parameters of a regression function  based on data for all 50 states in the United States or on data for all visits to a website. What is the interpretation of the estimated parameters and the standard errors? In practice, researchers typically assume that the sample is  randomly drawn  from a large population of interest and report standard errors that are designed to capture sampling variation. This is common even in applications where it is difficult to articulate what that population of interest is, and how it differs from the sample. In this article, we explore an alternative approach to inference, which is partly design-based. In a design-based setting, the values of some of the regressors can be manipulated, perhaps through a policy intervention. Design-based uncertainty emanates from lack of knowledge about the values that the regression outcome would have taken under alternative interventions. We derive standard errors that account for design-based uncertainty instead of, or in addition to, sampling-based uncertainty. We show that our  standard errors in general are smaller than the usual infinite-population sampling-based standard errors and provide conditions under which they coincide.
\end{abstract}

%\textbf{JEL Classification: C14, C21, C52}

%\textbf{Keywords:\ Regression Analyses, Standard Errors, Confidence Intervals, Random Sampling, Random Assignment, Finite Population, Large Samples, Potential Outcomes, Causality}

\end{titlepage}

\baselineskip=20pt\newpage \setcounter{page}{1}\renewcommand{\thepage}{%
\arabic{page}}\renewcommand{\theequation}{\arabic{section}.\arabic{equation}}

\section{Introduction}

\label{section:introduction}

In the dominant approach to inference in the social
sciences uncertainty about population parameters is induced by random sampling from the population.
Moreover, it is typically assumed that the sample comprises only a small fraction of the population of interest. This
perspective is a natural and attractive one in many instances. For example, if one analyzes
individual-level data from the U.S. Current Population Survey, the Panel
Study of Income Dynamics, or the one percent public-use sample from the U.S. Census,
it is natural to regard the sample as a small random subset of the
population of interest. In many other settings, however, this sampling
perspective is less appropriate. For example, \citet{manski2018right} write
``Random sampling assumptions, however, are not natural when considering
states or counties as units of observation.'' In this article, we provide an
alternative framework for the interpretation of uncertainty in regression
analysis regardless of whether a substantial fraction of the population or
even the entire population is included in the sample. While our framework
accommodates sampling-based uncertainty, it also takes into account
design-based uncertainty, which arises when the parameter of interest is
defined in terms of the unobserved outcomes that some units would attain
under a certain intervention. Design-based uncertainty is often explicitly
accounted for in the analysis of randomized experiments where it is the
basis of randomization inference (\citeauthor{neyman1923}, %
\citeyear{neyman1923}; \citeauthor{rosenbaum_book}, \citeyear{rosenbaum_book}%
; \citeauthor{imbens2015causal}, \citeyear{imbens2015causal}), but it is
rarely explicitly acknowledged in regression analyses or, more generally, in
observational studies (exceptions include \citeauthor{samii2012equivalencies}, %
\citeyear{samii2012equivalencies}; \citeauthor{Freedman2008}, %
\citeyear{Freedman2008}; %
\citeauthor{lu2016randomization}, \citeyear{lu2016randomization}; %
\citeauthor{lin}, \citeyear{lin}).

To illustrate the differences between sampling-based inference and
design-based inference, consider two simple examples. In the example of
Table \ref{tabel_sampling}, there is a finite population consisting of $n$
units with each unit characterized by a pair of variables, $Y_i$ and $Z_i$.
Consider an estimand that is function of the full set of
pairs $\{(Y_i,Z_i)\}^n_{i=1}$. Uncertainty about such an estimand arises when we observe the values $(Y_i,Z_i)$
only for a sample, that is, for a subset of the population. In Table \ref%
{tabel_sampling}, inclusion of unit $i$ in a sample is coded by the binary
variable $R_i\in\{0,1\}$. An estimator is a function of the observed data, $\{(R_i,R_iY_i,R_iZ_i)\}^n_{i=1}$.
Sampling-based inference uses information about the process
that determines the sampling indicators $R_1,\ldots, R_n$ to assess the variability of estimators across different samples. The second and third  sets of columns in Table \ref{tabel_sampling} depicts such alternative samples.
\begin{table}[ht]
\caption{\textsc{: Sampling-based Uncertainty ($\checkmark$ is observed, $?$
is missing)}}
\label{tabel_sampling}\vskip1cm
\par
\begin{center}
\begin{tabular}{c|cccccccccccccc}
& \multicolumn{3}{c}{Actual} &  & \multicolumn{3}{c}{Alternative} &  &
\multicolumn{3}{c}{Alternative} &  & $\hdots$ &  \\
Unit & \multicolumn{3}{c}{Sample} &  & \multicolumn{3}{c}{Sample I} &  &
\multicolumn{3}{c}{Sample II} &  & $\hdots$ &  \\
& $Y_i$ & $Z_i$ & $R_i$ & \hskip0.5cm \ \ \  & $Y_i$ & $Z_i$ & $R_i$ & \hskip%
0.5cm \ \ \  & $Y_i$ & $Z_i$ & $R_i$ & \hskip0.5cm & $\hdots$ & $$ \\ \hline
&  &  &  &  &  &  &  &  &  &  &  &  &  &  \\
$1$ & $\checkmark$ & $\checkmark$ & 1 &  & $?$ & $?$ & 0 &  & $?$ & $?$ & 0
&  & $\hdots$ &  \\
$2$ & $?$ & $?$ & 0 &  & $?$ & $?$ & 0 &  & $?$ & $?$ & 0 &  & $\hdots$ &
\\
$3$ & $?$ & $?$ & 0 &  & $\checkmark$ & $\checkmark$ & 1 &  & $\checkmark$ &
$\checkmark$ & 1 &  & $\hdots$ &  \\
$4$ & $?$ & $?$ & 0 &  & $\checkmark$ & $\checkmark$ & 1 &  & $?$ & $?$ & 0
&  & $\hdots$ &  \\
\vdots & \vdots & \vdots & \vdots &  & \vdots & \vdots & \vdots &  & \vdots
& \vdots & \vdots &  & $\hdots$ &  \\
$n$ & $\checkmark$ & $\checkmark$ & 1 &  & $?$ & $?$ & 0 &  & $?$ & $?$ & 0
&  & $\hdots$ &
\end{tabular}%
\end{center}
\end{table}
Table \ref{tabel_assignment} depicts a different scenario in which we observe for each
unit in the population the value of one of two potential outcome variables, either $%
Y^*_i(1)$ or $Y^*_i(0)$, but not both. The binary variable $X_i\in\{0,1\}$ indicates which potential outcome we
observe. Consider an estimand that is a function of the full set
of triples $\{(Y^*_i(1),Y^*_i(0),X_i)\}^n_{i=1}$. As before, an estimator is a function
of the observed data, the pairs $(X_i,Y_i)$, for $i=1,\ldots,n$, where $Y_i=Y_i^*(X_i)$ is the realized value. Design-based inference uses information about the
process that determines the assignments $X_1,\ldots, X_n$ to assess the variability of
estimators across different samples. The second and third  sets of columns in Table \ref{tabel_assignment} depicts such alternative samples.
\begin{table}[ht]
\caption{\textsc{: Design-based Uncertainty ($\checkmark$ is observed, $?$
is missing)}}
\label{tabel_assignment}\vskip1cm
\par
\begin{center}
\begin{tabular}{c|cccccccccccccc}
& \multicolumn{3}{c}{Actual} &  & \multicolumn{3}{c}{Alternative} &  &
\multicolumn{3}{c}{Alternative} &  & $\hdots$ &  \\
Unit & \multicolumn{3}{c}{Sample} &  & \multicolumn{3}{c}{Sample I} &  &
\multicolumn{3}{c}{Sample II} &  & $\hdots$ &  \\
& $Y_i^{*}(1)$ & $Y_i^{*}(0)$ & $X_i$ & \hskip0.5cm \ \ \  & $Y_i^{*}(1)$ & $%
Y_i^{*}(0)$ & $X_i$ & \hskip0.5cm \ \ \  & $Y_i^{*}(1)$ & $Y_i^{*}(0)$ & $%
X_i $ & \hskip0.5cm & $\hdots$ & $$ \\ \hline
&  &  &  &  &  &  &  &  &  &  &  &  &  &  \\
$1$ & $\checkmark$ & $?$ & 1 &  & $\checkmark$ & $?$ & 1 &  & ? & \checkmark
& 0 &  & $\hdots$ &  \\
$2$ & $?$ & $\checkmark$ & 0 &  & $?$ & $\checkmark$ & 0 &  & ? & \checkmark
& 0 &  & $\hdots$ &  \\
$3$ & $?$ & $\checkmark$ & 0 &  & $\checkmark$ & $?$ & 1 &  & \checkmark & ?
& 1 &  & $\hdots$ &  \\
$4$ & $?$ & $\checkmark$ & 0 &  & $?$ & $\checkmark$ & 0 &  & \checkmark & ?
& 1 &  & $\hdots$ &  \\
\vdots & \vdots & \vdots & \vdots &  & \vdots & \vdots & \vdots &  & \vdots
& \vdots & \vdots &  & $\hdots$ &  \\
$n$ & $\checkmark$ & $?$ & 1 &  & $?$ & $\checkmark$ & 0 &  & ? & \checkmark
& 0 &  & $\hdots$ &
\end{tabular}%
\end{center}
\end{table}

More generally, we can have missing data processes that combine features of
these two examples, with some units not included in the sample at all, and
with some of the variables not observed for the sampled units. Articulating both
the exact nature of the estimand of interest and the source of uncertainty
that makes an estimator stochastic is a crucial first step to valid
inference. For this purpose, it will be useful to distinguish between
descriptive estimands, where uncertainty stems solely from not observing all
units in the population of interest, and causal estimands, where the
uncertainty stems, at least partially, from unobservability of some of the
potential outcomes.

The main formal contribution of this article is to generalize the results
for the approximate variance for multiple linear regression estimators
associated with the work by \citet{eicker}, \citet{huber}, and %
\citet{white1980robust,white1980using,white1982maximum}, EHW from
hereon, in two ways. First, our framework allows for sampling from a finite population, whereas the EHW results assume random sampling from an infinite population.
Second, 
 our framework explicitly takes into account
design-based uncertainty.
Incorporating these generalizations requires developing a new framework for
regression analysis, nesting as special cases the \citet{neyman1923} analysis of randomized experiments with binary treatments, as well as the generalizations to randomized experiments with additional regressors in 
\citet{samii2012equivalencies}, \citet{Freedman2008}, 
and \citet{lin}. We show that in large
samples the widely used EHW robust standard errors are
conservative, and only correct in special cases. Moreover, we show that the
presence of attributes -- that is, characteristics of the units fixed in our
repeated sampling thought experiments -- can be exploited to improve on the
EHW variance estimator, and we propose variance estimators that
do so. Another advantage of the formal separation into sampling-based and
design-based uncertainty is that it allows us to clarify the distinction 
between the assumptions needed for internal and external validity (\citeauthor{shadishcookcampbell}, %
\citeyear{shadishcookcampbell}; \citeauthor{manski2013public}, %
\citeyear{manski2013public}; \citeauthor{deaton2010}, \citeyear{deaton2010})
in terms of these two sources of uncertainty.

Our results are relevant in empirical settings where
researchers have a random sample from a finite population and the ratio of
the sample size to the population size is sufficiently large so that the proposed
finite-population correction matters. Examples of such settings include large-scale experiments
\citep[see][]{muralidharan2017experimentation}, settings where the cost of data acquisition
motivates the use of random samples
\citep[see, e.g.,][]{keels2005fifteen}, as well as analyses based on public-use census samples, like the 2010 Integrated Public
Use Microdata Series (IPUMS) data (which is a 10 percent sample of the
U.S.~Census). More importantly in our view, our results are relevant in empirical settings
where it is not natural to think of the data as a random sample from a
well-defined population. Instead the researcher may have the entire
population, e.g., states or counties as in the
\citet{manski2018right} quote, or the set of all visits to a website, or the
researcher may have a convenience sample. In that case our design-based
approach to uncertainty provides a coherent interpretation for
sampling-based standard errors. It also provides methods that
exploit the presence of attributes to calculate improved (i.e., less conservative) standard errors.

\section{A Simple Example}

\label{section:cases}

In this section we set the stage for the problems discussed in the current
article by discussing least squares estimation in a simple example with a
single binary regressor. We make four points. First, we show how
design-based uncertainty affects the variance of regression estimators.
Second, we show that the standard Eicker-Huber-White (EHW) variance
estimator remains conservative when we take into account design-based
uncertainty. Third, we show that there is a simple finite-population
correction to the EHW variance estimator for descriptive
estimands but not for causal estimands. Fourth, we discuss the relation
between the two sources of uncertainty and the notions of internal and
external validity. Proofs of the results in this section are in a
supplementary appendix.

We focus on a setting with a finite population of size $n$. We sample $N$
units from this population, with $R_i\in\{0,1\}$ indicating whether a unit
was sampled ($R_i=1$) or not $(R_i=0$), so that $N=\sum_{i=1}^n R_i$. There
is a single binary regressor, $X_i\in\{0,1\}$, and $n_x$ (resp.\ $N_x$) is
the number of units in the population (resp.\ the sample) with $X_i=x$. 
Units could be U.S. states and the binary regressor $X_i$ could be an
indicator for a state regulation, say the state having a right-to-carry law (RTC),
as in \citet{manski2018right} and \citet{donohue2017right}.
We view the regressor $X_i$ not as a fixed attribute or characteristic of each
unit, but  as a cause or policy variable whose value could have been
different from the observed value. This generates missing data of the type
shown in Table \ref{tabel_assignment}, where only some of the states of the
world are observed, implying that there is design-based uncertainty.
Formally, using the Rubin causal model or potential outcome framework (%
\citeauthor{neyman1923}, \citeyear{neyman1923}; %
\citeauthor{rubin1974estimating}, \citeyear{rubin1974estimating}; %
\citeauthor{holland1986statistics}, \citeyear{holland1986statistics}; %
\citeauthor{imbens2015causal}, \citeyear{imbens2015causal}), we postulate
the existence of two potential outcomes for each unit, denoted by
$Y_i^{*}(1) $ and $Y_i^{*}(0)$. For the RTC example, $Y_i^{*}(1) $
and $Y_i^{*}(0)$ could be state-level crime rates
with and without RTC. The realized outcome is
\begin{equation*}
Y_i=Y_i^{*}(X_i)= \left\{
\begin{array}{ll}
Y_i^{*}(1) & \mathrm{if}\ X_i=1, \\
Y_i^{*}(0) & \mathrm{if}\ X_i=0,%
\end{array}
\right.
\end{equation*}
which is the observed state-level crime rate in the RTC example.

In our setting, potential outcomes
are viewed as non-stochastic attributes for unit $i$, irrespective of the
realized value of $X_i$. They, as well as the additional observed
attributes, $Z_i$ allow for in later sections, remain fixed in repeated sampling thought experiments,
whereas $R_i$ and $X_i$ are stochastic and, as a result, so are the realized
outcomes  in the sample, $Y_i$. In the current section, we abstract from the presence of
fixed observed attributes, which will play an important role in Section \ref%
{main}. Let ${\mathbf{Y}}$, ${\mathbf{Y}}^{*}(1)$, ${\mathbf{Y}}^{*}(0)$, ${%
\mathbf{R}}$, and ${\mathbf{X}}$ be the population $n$-vectors with $i$-th element
equal to $Y_i$, $Y_i^{*}(1)$, $Y_i^{*}(0)$, $R_i$, and $X_i$ respectively.
For sampled units (units with $R_i=1$) we observe $X_i$ and $Y_i$. For all units we observe $R_i$.

In general, estimands are functions of the full set of  values $({%
\mathbf{Y}}^{*}(1),{\mathbf{Y}}^{*}(0),{\mathbf{X}},{\mathbf{R}})$ for all units in the population, both those in the sample and those not in the sample. We
consider two types of estimands, descriptive and causal. If an estimand can
be written as a function of $({\mathbf{Y}},{\mathbf{X}})$, free of
dependence on ${\mathbf{R}}$ and on the potential outcomes beyond the
realized outcome, we label it a \textit{descriptive} estimand. Intuitively a
descriptive estimand is an estimand whose value would be known with
certainty if we observe  the realized values of all variables for all
units in the population. If an estimand cannot be written as a function of $%
({\mathbf{Y}},{\mathbf{X}},{\mathbf{R}})$ because it depends on the
potential outcomes ${\mathbf{Y}}^{*}(1)$ and ${\mathbf{Y}}^{*}(0)$, then we
label it a \textit{causal} estimand.\footnote{This does not partition the set of estimands. For example, there could be estimands that are functions of $({\mathbf{Y}},{\mathbf{X}},{\mathbf{R}})$, but not of  $({\mathbf{Y}},{\mathbf{X}})$, although it is difficult to think of interesting ones that are.}

We now consider in our binary regressor example three closely related
estimands, one descriptive and two causal:
\begin{equation*}
\theta^{\mathrm{descr}} =\theta^\mathrm{descr}({\mathbf{Y}},{\mathbf{X}}) =%
\frac{1}{n_{1}}\sum_{i=1}^nX_iY_i^{}- \frac{1}{n_0}\sum_{i=1}^n
(1-X_i)Y_i^{},
\end{equation*}
\begin{equation*}
\theta^{\mathrm{causal},\mathrm{sample}} =\theta^{\mathrm{causal},\mathrm{%
sample}} ({\mathbf{Y}}^{*}(1),{\mathbf{Y}}^{*}(0),{\mathbf{R}}) = \frac{1}{N}%
\sum_{i=1}^n R_i\, \Bigl( Y_i^{*}(1)-Y_i^{*}(0)\Bigr),
\end{equation*}
and
\begin{equation*}
\theta^\mathrm{causal}=\theta^\mathrm{causal}({\mathbf{Y}}^{*}(1),{\mathbf{Y}%
}^{*}(0)) = \frac{1}{n}\sum_{i=1}^n \Bigl( Y_i^{*}(1)-Y_i^{*}(0)\Bigr).
\end{equation*}
We focus on the properties of a particular estimator:
\begin{equation*}
\widehat\theta =\frac{1}{N_1}\sum_{i=1}^n R_iX_iY_i^{}- \frac{1}{N_0}%
\sum_{i=1}^n R_i(1-X_i)Y^{}_i.
\end{equation*}
This is also the least squares estimator of the coefficient on $X_i$ for
the regression in the sample  of $Y_i^{}$ on $X_i$ and a constant. There are two sources of
randomness in this estimator: a sampling component arising from the
randomness of $\boldsymbol{R}$ and a design component arising from the
randomness of $\boldsymbol{X}$. We refer to the uncertainty generated by the
randomness in the sampling component as \textit{sampling-based uncertainty},
and the uncertainty generated by the design component as \textit{%
design-based uncertainty}.

Next we consider the first two moments of $\widehat\theta$ under
the following two assumptions:
\begin{assumption}
\label{assumption:randomsampling1} \textsc{(Random Sampling / External
Validity)}
\begin{equation*}
\Pr\left(\boldsymbol{R}=\boldsymbol{r}\right) = 1\biggl/\left(%
\begin{array}{c}
n \\
N%
\end{array}%
\right),
\end{equation*}
for all $n$-vectors ${{\mathbf{r}}}$ with $\sum_{i=1}^n r_{i}=N$.
\end{assumption}
\begin{assumption}
\label{assumption:randomass1} \textsc{(Random Assignment / Internal Validity)%
}
\begin{equation*}
\Pr\left({{\mathbf{X}}}={{\mathbf{x}}}|{\mathbf{R}}\right) = 1\biggl/\left(%
\begin{array}{c}
n \\
n_1%
\end{array}%
\right),
\end{equation*}
for all $n$-vectors ${{\mathbf{x}}}$ with $\sum_{i=1}^n X_i=n_1$.
\end{assumption}

We start by studying the first moment of the estimator, conditional on $%
(N_1,N_0)$, and only for the cases where $N_1\geq 1$ and $N_0\geq 1$. We
leave this latter conditioning implicit in the notation throughout this
section. 
We also condition implicitly on the fixed potential outcomes ${\mathbf{Y}}^{*}(1)$ and ${\mathbf{Y}}^{*}(0)$.
Taking the expectation only over the random sampling, or taking
the expectation only over the random assignment, or over both, we find:
\begin{equation}  \label{equation:condexpx}
E\big[\widehat{\theta}\,|\,{\mathbf{X}},N_1,N_0\big]=\theta^{\mathrm{descr}},
\end{equation}
\begin{equation}  \label{equation:condexpr}
E\big[\widehat{\theta}\,|\,{\mathbf{R}},N_1,N_0\big]=\theta^{\mathrm{causal},%
\mathrm{sample}},
\end{equation}
\begin{equation*}
E\big[\widehat{\theta}\,|\,N_1,N_0\big] =E\big[\theta^{\mathrm{descr}%
}\,|\,N_1,N_0\big] =E\big[\theta^{\mathrm{causal},\mathrm{sample}%
}\,|\,N_1,N_0\big] =\theta^{\mathrm{causal}}.
\end{equation*}

Next, we look at the variance of the estimator, maintaining both the random
assignment and random sampling assumption.
Define the population variances
\begin{equation*}
S^2_x=\frac{1}{n-1}\sum_{i=1}^n \left(Y_i^{*}(x)-\frac{1}{n}\sum_{j=1}^n
Y_j^{*}(x)\right)^2,\hskip1cm {\rm  for}\ x=0,1,
\end{equation*}
and%
\begin{equation*}
S^2_\theta=\frac{1}{n-1}\sum_{i=1}^n \left(Y_i^{*}(1)-Y_i^{*}(0)-\frac{1}{n}%
\sum_{j=1}^n (Y_j^{*}(1) -Y_j^{*}(0))\right)^2.
\end{equation*}
We consider the variance of $%
\widehat\theta$, as well as two conditional versions of this variance. We define the ``sampling variance'' conditional on ${\mathbf{X}}$, so that only the
sampling uncertainty is taken into account. Analogously, we define the ``design variance''
conditional on ${\mathbf{R}}$, so that only the design uncertainty is taken
into account. To make the different variances interpretable, we look at the expected value of the variances, taking the expectation both over the assignment and the sampling.
\begin{align}  \label{equation:totalvariance}
V^\mathrm{total}(N_1,N_0,n_1,n_0)&=\mbox{var}\big(\widehat{\theta}\,|\,
N_1,N_0\big)=\frac{S^2_1}{N_1}+\frac{S^2_0}{N_0}-\frac{S^2_{\theta}}{n_0+n_1},
\end{align}
\begin{align*}
V^\mathrm{sampling}(N_1,N_0,n_1,n_0)&=E\left[\mbox{var}\big(\widehat{\theta}%
\,|\, {\mathbf{X}},N_1,N_0\big)\,\big|\, N_1,N_0\right] =\frac{S^2_1}{N_1}
\left( 1-\frac{N_1}{n_1}\right) +\frac{S^2_0}{N_0}\left(1-\frac{N_0}{n_0}%
\right),
\end{align*}
\begin{align*}
V^{\mathrm{design}}(N_1,N_0,n_1,n_0)&=E\left[\mbox{var}\big(\widehat{\theta}%
\,|\, {\mathbf{R}},N_1,N_0\big)\,|\, N_1,N_0\right]=\frac{S^2_1}{N_1}+\frac{%
S^2_0}{N_0}-\frac{S^2_{\theta}}{N_0+N_1}.
\end{align*}

\begin{remark}
\textsc{Neyman Variance}\newline
The variance $V^{\mathrm{total}}(N_1,N_0,n_1,n_0)$ is the one derived by %
\citet{neyman1923} for randomized experiments. \hfill$\square$
\end{remark}

\begin{remark}
\textsc{Causal versus Descriptive Estimands}\newline
In general the variances $V^\mathrm{sampling}(N_1,N_0,n_1,n_0)$ and $V^{%
\mathrm{design}}(N_1,N_0,n_1,n_0)$ cannot be ranked: the sampling variance
can be very close to zero if the sampling rate $(N_0+N_1)/(n_0+n_1)$ is close to one, but it
can also be larger than the design variance if the sampling rate is small
and the variance of the treatment effect is substantial.\hfill$\square$
\end{remark}

\begin{remark}
\textsc{Infinite Population Case}\newline
If $n_0,n_1\rightarrow\infty$, the total variance and the  sampling variance are equal:
\[
\lim_{n_0,n_1,\rightarrow\infty}
V^\mathrm{total}(N_1,N_0,n_1,n_0)
=
\lim_{n_0,n_1,\rightarrow\infty}V^\mathrm{sampling}(N_1,N_0,n_1,n_0)
=\frac{S^2_1}{N_1}+\frac{S^2_0}{N_0}.
\]
In this case taking the design-based uncertainty into account does not matter. This result will be seen to carry over to more general cases in Section \ref{main}.
\hfill$\square$
\end{remark}

\begin{remark}
\textsc{Finite Population Correction}\newline
Whether the estimand is $\theta^\mathrm{causal}$ or $\theta^\mathrm{descr}$,
ignoring the fact that the population is finite generally leads to an
overstatement of the variance on average because it ignores the fact that we observe a
non-negligible share of the population:
\begin{equation*}
V^\mathrm{total}(N_1,N_0,\infty,\infty)-V^\mathrm{total}(N_1,N_0,n_1,n_0)=%
\frac{S^2_{\theta}}{n_0+n_1}\geq 0,
\end{equation*}
\begin{equation*}
V^\mathrm{sampling}(N_1,N_0,\infty,\infty)-V^\mathrm{sampling}%
(N_1,N_0,n_1,n_0)= \frac{S^2_1}{n_1} +\frac{S^2_0}{n_0} \geq 0.
\end{equation*}
If the estimand is $\theta^{\mathrm{causal},\mathrm{sample}}$, however, the
population size is irrelevant because units in the population but not in the
sample do not matter for the estimand:
\begin{equation*}
V^\mathrm{design}(N_1,N_0,\infty,\infty)= V^\mathrm{design}%
(N_1,N_0,n_1,n_0).\hfill\square
\end{equation*}
\end{remark}

\begin{remark}
\textsc{Internal versus External Validity}\newline
Often researchers are concerned about both the internal and external
validity of estimands and estimators (\citeauthor{shadishcookcampbell}, %
\citeyear{shadishcookcampbell}; \citeauthor{manski2013public}, %
\citeyear{manski2013public}; \citeauthor{deaton2010}, \citeyear{deaton2010}%
). The distinction between sampling and design-based uncertainty allows us
to clarify these concerns. Internal validity bears on the question of whether $%
\widehat\theta$ is a good estimator for $\theta^{\mathrm{causal},\mathrm{%
sample}}$. This relies on random assignment of the treatment. Whether or not
the sampling is random is irrelevant for this question because $\theta^{%
\mathrm{causal},\mathrm{sampling}}$ conditions on which units were sampled.
External validity bears on the question of whether $E[\theta^\mathrm{causal,sample}]$ is equal to
 $\theta^\mathrm{causal}$. This relies on the random sampling
assumption and does not require that the assignment is random. However, for $%
\widehat\theta$ to be a good estimator of $\theta^\mathrm{causal}$, which
is often the most interesting estimand, we need both internal and external
validity, and thus both random assignment and random sampling. \hfill$%
\square $
\end{remark}

In this single binary regressor example the \textrm{EHW} variance estimator can
be written as
\begin{equation*}
\widehat V^\mathrm{ehw}=\frac{N_1-1}{N^2_1} \widehat S_1^2 + \frac{N_0-1}{%
N_0^2} \widehat S_0^2, \hskip0.5cm \mathrm{where}\ \ \widehat S_1^2=\frac{1}{%
N_1-1}\sum_{i=1}^n R_iX_i\left(Y_i-\frac{1}{N_1}\sum_{i=1}^n
R_iX_iY_i\right)^2,
\end{equation*}
and $\widehat S_0^2$ is defined analogously. Adjusting the degrees of
freedom, using the modification proposed in \citet{mackinnon1985some}
specialized to this binary regressor example, we obtain $\widetilde V^\mathrm{%
ehw}=\widehat S_1^2/N_1+\widehat S_0^2/N_0,$ which is identical to the variance estimator
proposed by \citet{neyman1923}, with the expectation of this modified \textrm{EHW} variance estimator $\widetilde V^\mathrm{ehw}$  (conditional on $N_0$ and $%
N_1$) equal to the sampling variance in the infinite population case, $V^%
\mathrm{sampling}(N_1,N_0,\infty,\infty)$.  In the infinite population case the design-based uncertainty does not matter, so the EHW variance can be intepreted as implicitly taking into account design-uncertainty by focusing on the infinite population case.

We could also estimate the variance using resampling methods, which would give us variance estimates close to $\widehat{V}^{\mathrm{ehw}}$. To be precise, suppose we use the bootstrap where we
 draw $N_1$ bootstrap observations from the $N_1$ treated units
and $N_0$ bootstrap units from the $N_0$ control units. In that case the bootstrap variance would in expectation (over the bootstrap replications)
be equal to $\widehat{V}^{\mathrm{ehw}}$.

\begin{remark}
\textsc{Can We Improve on the EHW variance Estimator?}\newline
The difference between $E[\widetilde V^{\mathrm{ehw}}|N_1,N_0]$(or the Neyman variance)  and the total variance is equal to $S^2_\theta/n$.
The term $S^2_\theta$ is difficult to estimate because it depends on the unobserved differences $Y_i^{*}(1)-Y_i^{*}(0)$.
As a result,  $S^2_\theta/n$ is typically ignored in analyses of
randomized experiments \citep[see][]{imbens2015causal}. In particular, the EHW variance estimator
implicitly sets the estimator of $S^2_\theta$ to be equal to zero, resulting in conservative inference.
For the case of a randomized experiment with a binary treatment \cite{aronow2014sharp} provide  a lower bound for $S^2_\theta$ based on the Fr\'echet-Hoefffding inequality.
 In Section \ref{main},
we propose an improved variance estimator that exploits the presence of fixed
attributes.\hfill$\square$
\end{remark}

The appendix contains a Bayesian version of the analysis ofthe example from  this section. Similar to the results in this section we show that, when the estimand of interest is defined for a finite population, the posterior variance depends not only on the sample sizes for treated and non-treated but also on the respective population sizes. Also similar to the analysis of this section, the posterior variance formula depends on whether the estimand is descriptive or causal.

\section{The General Case}
\label{main}

This section contains the main formal results in the article. We focus on a
regression setting where we estimate a linear regression function for a
scalar outcome and a number of regressors. The setting we consider here
allows for the presence of two types of regressors. First, regressors that
are causal, in the sense that they generate potential outcomes. Second,
regressors that are attributes, in the sense that they are kept fixed for
each unit in the thought experiment that provides the basis for inference.
Which regressors are viewed as causal and which are viewed as attributes
depends on the interpretation we wish to give to the
regression estimates. If we wish to give a coefficient a causal
interpretation, the corresponding regressor must be a cause. If a regressor
is an attribute, the corresponding coefficient is simply estimating a
population difference between subpopulations of units. For example, if we
regress earnings on years of education, years of education may be the causal
variable of interest. On the other hand, if we regress earnings on an
indicator for participation in a job search program, age, and years of
education, then the indicator for the program participation may be viewed as the causal variable of
interest and age and years of education may be viewed as  attributes.
Because the repeated sampling thought experiment treats causes different
from attributes, the variance of the regression estimator will depend on
this designation.

\subsection{Set Up}

Consider a sequence of finite populations indexed by population size, $n$.
Unit $i$ in population $n$ is characterized by a set of fixed attributes $%
Z_{n,i}$ (including an intercept) and by a potential outcome function, $%
Y_{n,i}^{*}(\cdot)$, which maps causes, $U_{n,i}$, into outcomes, with the realized outcome denoted by  $%
Y_{n,i}^{}=Y_{n,i}^{*}(U_{n,i})$. $Z_{n,i}$ and $U_{n,i}$ are real-valued
column vectors, and $Y_{n,i}^{}$ is scalar. We do not place restrictions
on the types of the variables: they can be continuous, discrete, or mixed.
The attributes and potential outcome functions remain fixed in our repeated
sampling thought experiments. The realized outcomes in the sample vary from
sample to sample because the units in the sample and the values of the causal variables change.

There is a sequence of samples associated with the sequence of populations. We
will use $R_{n,i}=1$ to indicate that unit $i$ of population $n$ is sampled,
and $R_{n,i}=0$ to indicate that it is not sampled. For each unit in sample $%
n$, we observe the triple, $(Y_{n,i},U_{n,i},Z_{n,i})$.
Relative to Section \ref%
{section:cases} we now allow for more complicated assignment
mechanisms. In particular, we relax the assumption that the causes have
identical distributions.

\begin{assumption}
\label{assumption:assignment}\textsc{(Assignment Mechanism)} The assignments
$U_{n,1},\ldots ,U_{n,n}$ are jointly independent, and independent of $%
R_{n,1},\ldots ,R_{n,n}$, but not (necessarily) identically distributed
(i.n.i.d.).
\end{assumption}
This assumption assumes independence of the treatment assignments. This
is somewhat in contrast to the example in Section \ref{section:cases}, where
we fixed the marginal distribution of the regressor, allowing us to obtain exact finite sample results. We do not need this
here because we are focused on asymptotic results.
We can allow for some dependence in the assignment mechanism, {\it e.g.,} clustering of the type analyzed in \cite{abadie2017should}.

For what follows, it is convenient to work with a transformation $%
X_{n,1},\ldots ,X_{n,n}$ of $U_{n,1},\ldots ,U_{n,n}$ that removes the
correlation with the attributes.
%:\begin{equation}  \label{equation:uncorrelated}E\left[\sum_{i=1}^n X_{n,i}Z_{n,i}^{\prime }\right]=\sum_{i=1}^nE[X_{n,i}]Z_{n,i}^{\prime }=0,\end{equation}where the expectations are taken over the assignments.
This can be accomplished in the following way. We assume that the population
matrix $\sum_{i=1}^n Z_{n,i}Z_{n,i}^{\prime }$ is full-rank.
Then, define %equation (\ref{equation:uncorrelated}) holds for
\begin{equation}  \label{equation:transformation}
X_{n,i}=U_{n,i}-\Lambda_n Z_{n,i} \hskip1cm \mathrm{where}\ \ \Lambda_n =
\left(\sum_{i=1}^n E[U_{n,i}]Z_{n,i}^{\prime }\right)\left(\sum_{i=1}^n
Z_{n,i}Z_{n,i}^{\prime }\right)^{-1}.
\end{equation}
Later we formally make an assumption that will guarantee that this transformation is well-defined for large $n$. It is important to notice that, because $\Lambda_n
Z_{n,i}$ is deterministic in our setting and $U_{n,1}, \ldots, U_{n,n}$ are
i.n.i.d., the variables $X_{n,1}, \ldots, X_{n,n}$ are i.n.i.d.\ too.

For population $n$, let $\boldsymbol{Y}_n$, $\boldsymbol{X}_n$, $\boldsymbol{%
Z}_n$, $\boldsymbol{R}_n$, and $\boldsymbol{Y}^*_n(\cdot)$ be matrices that
collect outcomes, causes, attributes, sampling indicators, and potential
outcome functions.
We analyze the properties of the estimator $\widehat\theta_n$
obtained by minimizing least square errors in the sample:
\begin{equation}
\label{equation:least_squares}
(\widehat\theta_n,\widehat\gamma_n)= {\arg\min}_{(\theta,\gamma)}\,
\sum_{i=1}^n R_{n,i}\big(Y_{n,i}-X_{n,i}^{\prime }\theta-Z_{n,i}^{\prime
}\gamma\big)^2.
\end{equation}
The properties of the population regression residuals, $%
e_{n,i}=Y_{n,i}-X_{n,i}^{\prime }\theta_n-Z_{n,i}^{\prime }\gamma_n$, depend
on the exact nature of the estimands, $(\theta_n,\gamma_n)$. In what
follows, we will consider alternative target parameters, which in turn will
imply different properties for $e_{n,i}$. Notice also that, although the
transformation in (\ref{equation:transformation}) is typically unfeasible
(because the values of $E[U_{n,i}]$ may not be known), $\widehat\theta_n$ is
not affected by the transformation in the sense that the least squares
estimators $(\widetilde\theta_n,\widetilde\gamma_n)$, defined as
\begin{equation*}
(\widetilde\theta_n,\widetilde\gamma_n)= {\arg\min}_{(\theta,\gamma)}\,
\sum_{i=1}^n R_{n,i}\big(Y_{n,i}-U_{n,i}^{\prime }\theta-Z_{n,i}^{\prime
}\gamma\big)^2,
\end{equation*}
satisfy $\widehat\theta_n=\widetilde\theta_n$ (although, in general, $%
\widehat\gamma_n\neq \widetilde\gamma_n$). As a result, we can analyze the
properties of $\widehat\theta_n$, focusing on the properties of the
regression on $X_{n,1},\ldots, X_{n,n}$ instead of on $U_{n,1},\ldots,
U_{n,n}$.

We assume random sampling, with some conditions on the sampling rate to
ensure that the sample size increases with the population size.

\begin{assumption}
\label{assumption:randomsampling} \textsc{(Random Sampling)} (i) There is a
sequence of sampling probabilities, $\rho_n$, such that
\begin{equation*}
\Pr\left(\boldsymbol{R}_n=\boldsymbol{r}\right) =\rho_n^{\,\sum_{i=1}^n r_i}
\left(1-\rho_n\right) ^{n-\sum_{i=1}^n r_i}
\end{equation*}
for all $n$-vectors $\boldsymbol{r}$ with $i$-th element $r_i\in \{0,1\}$.
(ii) The sequence of sampling rates, $\rho_n$, satisfies $%
n\rho_n\rightarrow\infty$ and $\rho_n\rightarrow\rho\in [0,1]$.
\end{assumption}
The first part of
Assumption \ref{assumption:randomsampling}$(ii)$ guarantees that as the
population size increases, the (expected) sample size also increases. The
second part of Assumption \ref{assumption:randomsampling}$(ii)$ allows for
the possibility that, as $n$ increases, the sample size becomes a negligible fraction
of the population size so that the \textrm{EHW} results, corresponding to $\rho=0$, are included as a special case of our results.

The next assumption is a regularity condition bounding moments.

\begin{assumption}
\textsc{(Moments)} \label{assumption:moments} There exists some $\delta>0$
such that the sequences
\begin{equation*}
\frac{1}{n}\sum_{i=1}^n E[|Y_{n,i}|^{4+\delta}],\qquad \frac{1}{n}%
\sum_{i=1}^n E[\|X_{n,i}\|^{4+\delta}],\qquad\frac{1}{n}\sum_{i=1}^n
\|Z_{n,i}\|^{4+\delta}
\end{equation*}
are uniformly bounded.
\end{assumption}

Let
\begin{equation*}
W_n = \frac{1}{n}\sum_{i=1}^n \left(%
\begin{array}{c}
Y_{n,i} \\
X_{n,i} \\
Z_{n,i}%
\end{array}%
\right) \left(%
\begin{array}{c}
Y_{n,i} \\
X_{n,i} \\
Z_{n,i}%
\end{array}%
\right)^{\prime },\quad \Omega_n=\frac{1}{n}\sum_{i=1}^n E\left[\left(%
\begin{array}{c}
Y_{n,i} \\
X_{n,i} \\
Z_{n,i}%
\end{array}%
\right) \left(%
\begin{array}{c}
Y_{n,i} \\
X_{n,i} \\
Z_{n,i}%
\end{array}%
\right)^{\prime }\right].
\end{equation*}
So $\Omega_n=E[W_n]$, where the expectation is taken over the distribution
of $\boldsymbol{X}_n$. We also consider sample counterparts of $W_n$
and $\Omega_n$:
\begin{equation*}
\widetilde W_n = \frac{1}{N}\sum_{i=1}^N R_{n,i}\left(%
\begin{array}{c}
Y_{n,i} \\
X_{n,i} \\
Z_{n,i}%
\end{array}%
\right) \left(%
\begin{array}{c}
Y_{n,i} \\
X_{n,i} \\
Z_{n,i}%
\end{array}%
\right)^{\prime },\quad \widetilde\Omega_n = \frac{1}{N}\sum_{i=1}^n R_{n,i}
E\left[\left(%
\begin{array}{c}
Y_{n,i} \\
X_{n,i} \\
Z_{n,i}%
\end{array}%
\right) \left(%
\begin{array}{c}
Y_{n,i} \\
X_{n,i} \\
Z_{n,i}%
\end{array}%
\right)^{\prime }\right],
\end{equation*}
where $\widetilde\Omega_n =E[\widetilde W_n|{\mathbf{R}}_n]$. We will use
superscripts to indicate submatrices. For example,
\begin{equation*}
W_n=\left(%
\begin{array}{ccc}
W_n^{YY} & W_n^{YX} & W_n^{YZ} \\
W_n^{XY} & W_n^{XX} & W_n^{XZ} \\
W_n^{ZY} & W_n^{ZX} & W_n^{ZZ}%
\end{array}%
\right),
\end{equation*}
with analogous partitions for $\Omega_n$, $\widetilde W_n$, and $%
\widetilde\Omega_n$. Notice that the transformation in (\ref%
{equation:transformation}) implies that $\Omega_n^{XZ}$ and $\Omega_n^{ZX}$
are matrices will all zero entries.

We first obtain convergence results for the sample objects, $\widetilde W_n$
and $\widetilde \Omega_n$.

\begin{lemma}
\label{lemma:convergence2moments} Suppose Assumptions \ref%
{assumption:assignment}-\ref{assumption:moments} hold. Then, $\widetilde W_n
- \Omega_n \overset{p}{\rightarrow} 0$ , $\widetilde \Omega_n - \Omega_n
\overset{p}{\rightarrow} 0$ and $\widetilde W_n-W_n\overset{p}{\rightarrow}
0 $.
\end{lemma}
\noindent See appendix for proofs.

The next assumption imposes (deterministic) convergence of the expected value of the
second moments in the population.

\begin{assumption}
\textsc{(Convergence of moments)} \label{assumption:convergence} $%
\Omega_n\rightarrow \Omega$, which is full rank.
\end{assumption}

\subsection{Descriptive and Causal Estimands}

\label{subsection:descriptiveandcausal}

We now define the descriptive and causal estimands that generalize $\theta ^{%
\mathrm{descr}}$, $\theta ^{\mathrm{causal},\mathrm{sample}}$, and $\theta^{%
\mathrm{causal}}$ from Section \ref{section:cases} to a regression context.

\begin{definition}
\textsc{Causal and Descriptive Estimands}\newline
For a given population $n$, with potential outcome functions $\boldsymbol{Y}%
_n^{*}(\cdot)$, causes $\boldsymbol{X}_n$, attributes $\boldsymbol{Z}_n$,
and sampling indicators $\boldsymbol{R}_n$:\vspace*{-5pt}

\begin{itemize}
\setlength\itemsep{-2pt}

\item[(i)] Estimands are functionals of $(\boldsymbol{Y}_n^{*}(\cdot),%
\boldsymbol{X}_n ,\boldsymbol{Z}_n,\boldsymbol{R}_n)$, permutation-invariant
in the rows of the arguments.

\item[(ii)] Descriptive estimands are estimands that can be written in terms
of $\boldsymbol{Y}_n^{}$, $\boldsymbol{X}_n$, and $\boldsymbol{Z}_n$, free
of dependence on $\boldsymbol{R}_n$, and free of dependence on $\boldsymbol{Y%
}_n^{*}(\cdot)$ beyond dependence on $\boldsymbol{Y}_n^{}$.

\item[(iii)] Causal estimands are estimands that cannot be written in terms
of $\boldsymbol{Y}_n^{}$, $\boldsymbol{X}_n$, $\boldsymbol{Z}_n$, and $%
\boldsymbol{R}_n$, because they depend on the potential outcome functions $%
\boldsymbol{Y}_n^*(\cdot)$ beyond the realized outcomes, $\boldsymbol{Y}_n$.
\end{itemize}
\end{definition}

Causal estimands depend on the values of potential outcomes beyond the
values that can be inferred from the realized outcomes. Given a sample, the
only reason we may not be able to infer the exact value of a descriptive estimand
is that we do not see all the units in the population. In contrast, even if
we observe all units in a population, we are unable to infer the value of a
causal estimand because its value depends on potential outcomes.

We define three estimands of interest, which under the conditions above
exist with probability approaching one:
\begin{equation}  \label{equation:descriptive}
\left(%
\begin{array}{c}
\theta^{\mathrm{descr}}_n \\
\gamma^{\mathrm{descr}}_n%
\end{array}%
\right)=\left(%
\begin{array}{cc}
W^{XX}_n & W^{XZ}_n \\
W^{ZX}_n & W^{ZZ}_n%
\end{array}%
\right)^{-1}\left(%
\begin{array}{c}
W^{XY}_n \\
W^{ZY}_n%
\end{array}%
\right),
\end{equation}
\begin{equation}  \label{equation:causalsample}
\left(%
\begin{array}{c}
\theta^{\mathrm{causal}, \mathrm{sample}}_n \\
\gamma^{\mathrm{causal}, \mathrm{sample}}_n%
\end{array}%
\right)=\left(%
\begin{array}{cc}
\widetilde \Omega^{XX}_n & \widetilde\Omega^{XZ}_n \\
\widetilde\Omega^{ZX}_n & \widetilde\Omega^{ZZ}_n%
\end{array}%
\right)^{-1}\left(%
\begin{array}{c}
\widetilde\Omega^{XY}_n \\
\widetilde\Omega^{ZY}_n%
\end{array}%
\right),
\end{equation}
and
\begin{equation}  \label{equation:causal}
\left(%
\begin{array}{c}
\theta^{\mathrm{causal}}_n \\
\gamma^{\mathrm{causal}}_n%
\end{array}%
\right)=\left(%
\begin{array}{cc}
\Omega^{XX}_n & \Omega^{XZ}_n \\
\Omega^{ZX}_n & \Omega^{ZZ}_n%
\end{array}%
\right)^{-1}\left(%
\begin{array}{c}
\Omega^{XY}_n \\
\Omega^{ZY}_n%
\end{array}%
\right).
\end{equation}
Alternatively, the estimands in (\ref{equation:descriptive}) to (\ref%
{equation:causal}) can be defined as the coefficients
that correspond to the orthogonality conditions in terms of the residuals $%
e_{n,i}=Y_{n,i}-X_{n,i}^{\prime }\theta_n-Z^{\prime }_{n,i}\gamma_n$,
\begin{equation*}
\frac{1}{n}\sum_{i=1}^n \left(%
\begin{array}{c}
X_{n,i} \\
Z_{n,i}%
\end{array}%
\right)e_{n,i}=0,\quad \frac{1}{n}\sum_{i=1}^n R_{n,i}E\left[\left(%
\begin{array}{c}
X_{n,i} \\
Z_{n,i}%
\end{array}%
\right)e_{n,i}\right]=0,\quad \frac{1}{n}\sum_{i=1}^n E\left[\left(%
\begin{array}{c}
X_{n,i} \\
Z_{n,i}%
\end{array}%
\right)e_{n,i}\right]=0,
\end{equation*}
for the descriptive, causal-sample, and causal estimands respectively. We
will study the properties of the least squares estimator, $\widehat\theta_n$%
, defined by
\begin{equation*}
\left(%
\begin{array}{c}
\widehat\theta_n \\
\widehat\gamma_n%
\end{array}%
\right)=\left(%
\begin{array}{cc}
\widetilde W^{XX}_n & \widetilde W^{XZ}_n \\
\widetilde W^{ZX}_n & \widetilde W^{ZZ}_n%
\end{array}%
\right)^{-1}\left(%
\begin{array}{c}
\widetilde W^{XY}_n \\
\widetilde W^{ZY}_n%
\end{array}%
\right),
\end{equation*}
as an estimator of the parameters defined in equations (\ref%
{equation:descriptive}) to (\ref{equation:causal}).

Notice that
%, by the law of total expectation and because potential outcome functions are fixed in our framework, 
$\theta^{\mathrm{causal}, \mathrm{sample}}_n$
and $\theta^{\mathrm{causal}}_n$ are causal estimands, while $\theta^{\mathrm{descr}%
}_n$ is not. However, the fact that an estimand is
causal according to our definition does not imply it has an interpretation
as an average causal effect. In Section \ref{section:causal} we present
conditions under which the regression estimand does have such an
interpretation.

\subsection{Causal Interpretations of the Estimands}

\label{section:causal}

By construction, the descriptive estimand can be interpreted as the set of
coefficients of a population best linear predictor (least squares),  e.g.,
\citeauthor{goldberger1991course}
(\citeyear{goldberger1991course}). A more
challenging question concerns the interpretation of the two causal
estimands, and in particular their relation to the potential outcome
functions. In this section we investigate this question.

The next assumption generalizes random assignment to allow for some dependence of the assignment $U_{n,i}$ on the attributes $Z_{n,i}$.

\begin{assumption}
\label{assumption:unconf}\textsc{(Expected Assignment)}
$(i)$
There exists a sequence of functions $h_n$ such that
\begin{equation*}
E[U_{n,i}]=h_n(Z_{n,i}),
\end{equation*}
and $(ii)$ there exists a sequence of matrices $B_n$ such that for all $z$
\begin{equation*}
h_n(z)=B_n z,
\end{equation*}
for all $n$ large enough.
\end{assumption}

Assumption \ref{assumption:unconf} looks very different from conventional
exogeneity or unconfoundedness conditions, where the residuals are assumed
to be (mean-) independent of the regressors, and so it merits some
discussion. First, note that if the treatment $U_{n,i}$ is
randomly assigned $E[U_{n,i}]$ is constant and Assumption \ref%
{assumption:unconf}$(i)$ and $(ii)$ are automatically satisfied as long as $Z_{n,i}$ includes
an intercept.

 Formally, Assumption \ref{assumption:unconf} relaxes the
completely randomized assignment setting by allowing the distribution of $%
U_{n,i}$ to depend on the attributes. However, this dependence is restricted
in that the mean of $U_{n,i}$ is linear in $Z_{n,i}$. For example,
Assumption \ref{assumption:unconf} holds automatically when $U_{n,1},\ldots,
U_{n,n}$ are identically distributed and $Z_{n,i}$ contains a saturated set
of indicators for all possible values of the attributes.

In the special case where the treatment is binary and $E[U_{n,i}]$ is the
propensity score, the assumption amounts to combination of an
unconfoundedness assumption that the treatment assignment does not depend on
the potential outcomes and a linear model for the
propensity score.

Later in this section, we will show that under a set of conditions that
includes Assumption \ref{assumption:unconf}, the two estimands $\theta_n^\mathrm{causal}$ and $%
\theta^{\mathrm{causal},\mathrm{sample}}$ can be interpreted as weighted
averages of unit-level causal effects. The connection between linearity in
the propensity score (\citet{rosenbaum1983central}), in our analysis
represented by $E[U_{n,i}]=B_n Z_{n,i}$, and the interpretation of
population regression coefficients as weighted averages of heterogeneous
causal effects has been previously noticed in related contexts (see %
\citeauthor{angrist1998}, \citeyear{angrist1998}; \citeauthor{angristpischke}%
, \citeyear{angristpischke}; \citeauthor{aronow2016}, \citeyear{aronow2016}; %
\citeauthor{tymon2017}, \citeyear{tymon2017}).

\begin{assumption}
\label{assumption:linear}\textsc{(Linearity of Potential Outcomes)} For all $%
u$,
\begin{equation}
\label{equation:pot_outcomes}
Y_{n,i}^*(u)=u^{\prime }\theta_{n,i}+\xi_{n,i},
\end{equation}
where $\theta_{n,i}$ and $\xi_{n,i}$ are non-stochastic.
\end{assumption}

In this formulation, any dependence of the potential outcomes $%
Y_{n,i}^{*}(u) $ on observed or unobserved attributes is subsumed by $%
\theta_{n,i}$ and $\xi_{n,i}$, which are non-stochastic. Each element of the
vector $\theta_{n,i}$ represents the causal effect of
increasing the corresponding value of $U_{n,i}$ by one unit.

The linearity in Assumption \ref{assumption:linear} is a strong restriction
in many settings. However, in some leading cases -- in particular, when the causal
variable is binary or, more generally when the causal variable takes on only
a finite number of values -- one can ensure that this assumption holds by
including in $U_{n,i}$ indicator variables representing each but one of the
possible values of the cause. With Assumption \ref{assumption:linear}
we are able to provide a
more transparent interpretation of the regression estimator. 
%Without it, but maintaining the other assumptions, we can only interpret the regression estimator as a weighted average derivative.

\begin{theorem}
\label{theorem:causal} Suppose Assumptions \ref{assumption:assignment}-\ref%
{assumption:linear} hold. Then, for all $n$ large enough,
\begin{align*}
\theta_n^{\mathrm{causal}}&=\left(\sum_{i=1}^n E\left[W_{n,i}^{XX} \right]%
\right)^{-1}\sum_{i=1}^nE\left[W_{n,i}^{XX}\right]\theta_{n,i}, \\
\shortintertext{and, with probability approaching one,} \theta_n^{\mathrm{%
causal}, \mathrm{sample}}&=\left(\sum_{i=1}^n R_{n,i}E\left[W_{n,i}^{XX} %
\right]\right)^{-1}\sum_{i=1}^n R_{n,i}E\left[W_{n,i}^{XX}\right]%
\theta_{n,i},
\end{align*}
where $W_{n,i}^{XX}=X_{n,i}X_{n,i}^{\prime }$.
\end{theorem}

\begin{theorem}
\label{theorem:nonlinear} Suppose that Assumptions \ref%
{assumption:assignment}-\ref{assumption:unconf} hold. Moreover, assume that $%
X_{n,1}, \ldots, X_{n,n}$ are continuous random variables with convex and
compact supports, and that the potential outcome functions, $Y^*_{n,i}(\cdot)$ are
continuously differentiable. Then, there exist random variables $%
v_{n,1},\ldots, v_{n,n}$ such that, for $n$ sufficiently large,
\begin{align*}
\theta_n^{\mathrm{causal}}&=\left(\sum_{i=1}^n E\left[W_{n,i}^{XX} \right]%
\right)^{-1}\sum_{i=1}^nE\left[W_{n,i}^{XX}\varphi_{n,i}\right], \\
\shortintertext{and, with probability approaching one,} \theta_n^{\mathrm{causal}, \mathrm{sample}%
}&=\left(\sum_{i=1}^n R_{n,i}E\left[W_{n,i}^{XX} \right]\right)^{-1}%
\sum_{i=1}^n R_{n,i}E\left[W_{n,i}^{XX}\varphi_{n,i}\right],
\end{align*}
where $\varphi_{n,i}$ is the derivative of $Y_{n,i}^{*}(\cdot)$ evaluated at
$v_{n,i}$.
\end{theorem}

\begin{remark}
Here, we provide a simple example that shows how the result in Theorems \ref%
{theorem:causal} and \ref{theorem:nonlinear} may not hold in the absence of
Assumption \ref{assumption:unconf}. Consider the population with three units
described in Table \ref{table:example} . For simplicity, we drop the
subscript $n$. We can also add replicates of these observations to make the example hold for any population size. In this example, $E[U_{i}]=3b Z_i^2-2b$ is a non-linear
function of $Z_{i}$. Notice that
\begin{equation*}
\sum_{i=1}^3 E[U_i]/3=\sum_{i=1}^3 E[U_i] Z_i/3=0,
\end{equation*}
so that $X_i=U_i$. Therefore, $E[X_i^2]=E[U_i^2]$. Also, because potential
outcomes do not depend on $X_i$, it follows that $%
E[X_iY_i^{}]=E[X_i]Y^*_i(1)=E[U_i]Y^*_i(1)$. As a result,
\begin{equation*}
\theta^{\mathrm{causal}}= \left(\sum_{i=1}^3
E[X_i^2]\right)^{-1}\sum_{i=1}^3 E[X_i]Y^*_i(1)=\frac{ab}{2b^2+1},
\end{equation*}
which is different from zero as long as $ab\neq 0$. In this example, all the
potential outcome functions $Y_i^{*}(\cdot)$ are flat as a function of $x$,
so all unit-level causal effects of the type $Y_i^{*}(u)-Y_i^{*}(u^{\prime
}) $ are zero, and yet the causal least squares estimand can be positive or
negative depending on the values of $a$ and $b$.\hfill$\square$
\end{remark}\vspace*{-0.3cm}
\begin{table}[ht]
\caption{: An Artificial Example}
\label{table:example}\vskip1cm
\par
\begin{center}
\begin{tabular}{c|rrrc}
Unit & $Y_{i}^{*}(u)$ & $Z_i$ & $E[U_i]$ & $\mbox{var}(U_i)$ \\ \hline
&  &  &  &  \\[-2ex]
$1$ & $a$ & $-1$ & $b$ & $1$ \\
$2$ & $0$ & $0$ & $-2b$ & $1$ \\
$3$ & $2a$ & $1$ & $b$ & $1$%
\end{tabular}%
\end{center}
\end{table}\vspace*{-.8cm}
\subsection{The Asymptotic Distribution of The Least Squares Estimator}
\label{section:distribution}

In this section we present the main result of the article,
describing the properties of the least squares estimator viewed as an
estimator of the causal estimands and, separately, viewed as an estimator of
the descriptive estimand. In contrast to Section \ref{section:cases}, we do
not have exact results, relying instead on asymptotic results based on
sequences of populations.

First, we define the population residuals, denoted by $\varepsilon _{n,i}$,
relative to the population causal estimands,
\begin{equation}  \label{equation:causalres}
\varepsilon _{n,i}=Y_{n,i}-X_{n,i}^{\prime }\theta _n^\mathrm{causal}%
-Z_{n,i}^{\prime }\gamma _n^\mathrm{causal}.
\end{equation}

\begin{remark}
The definition of the residuals, $\varepsilon_{n,1}, \ldots,
\varepsilon_{n,n}$, mirrors that in conventional regression analysis, but
their properties are conceptually different. For instance, the residuals
need not be stochastic. If they are stochastic, they are so because of their
dependence on $\boldsymbol{X}_n$.\hfill$\square$
\end{remark}

\begin{remark} We define the residuals here with respect to the population causal parameters $\theta _n^\mathrm{causal}$ and
$\gamma _n^\mathrm{causal}$. Because we focus here on asymptotic results, the difference between the causal and descriptive parameters vanishes, and so
defining the residuals in terms of the descriptive parameters would lead to the same results.\hfill$\square$
\end{remark}

Under the assumption that the $X_{n,i}$ are jointly independent (but not
necessarily identically distributed), the $n$ products $X_{n,i} \varepsilon
_{n,i}$ are also jointly independent but not identically distributed. Most
importantly, in general the expectations $E[X_{n,i} \varepsilon _{n,i}] $
may vary across $i$, and need not all be zero. However, as shown in Section %
\ref{subsection:descriptiveandcausal}, the \textit{averages} of these
expectations over the entire population are guaranteed to be zero by the
definition of $(\theta^\mathrm{causal}_n,\gamma^\mathrm{causal}_n)$. Define
the limits of the population variance,
\begin{equation*}
\Delta^{\mathrm{cond}}=\lim_{n\rightarrow \infty}\frac{1}{n}\sum_{i=1}^{n}%
\mbox{var}\left(X_{n,i} \varepsilon _{n,i}\right),
\end{equation*}
and the expected outer product
\begin{equation*}
\Delta^{\mathrm{ehw}}=\lim_{n\rightarrow \infty}\frac{1}{n}\sum_{i=1}^{n} E%
\left[\varepsilon_{n,i}^2X_{n,i}X_{n,i}^{\prime }\right].
\end{equation*}
The difference between $\Delta^\mathrm{ehw}$ and $\Delta^\mathrm{cond}$ is
the limit of the average outer product of the means,
\begin{equation*}
\Delta^\mu=\Delta^\mathrm{ehw}-\Delta^\mathrm{cond}=\lim_{n\rightarrow
\infty}\frac{1}{n}\sum_{i=1}^{n}
E[X_{n,i}\varepsilon_{n,i}]E[X_{n,i}\varepsilon_{n,i}]^{\prime },
\end{equation*}
which is positive semidefinite. We assume existence of these limits.

\begin{assumption}
\textsc{(Existence of Limits)} \label{assumption:limits} $\Delta^\mathrm{cond%
}$ and $\Delta^\mathrm{ehw}$ exist and are positive definite.
\end{assumption}

\begin{theorem}
\label{theorem:asympdist} Suppose Assumptions \ref{assumption:assignment}-%
\ref{assumption:limits} hold, and let $\Gamma  =
\Omega^{XX}=\lim_{n\rightarrow\infty} \Omega^{XX}_n$. Then, \vspace*{-10pt}

\begin{itemize}
\setlength\itemsep{-2pt}

\item[(i)]
\begin{equation*}
\sqrt{N}\left( \widehat{\theta}_n-\theta_n^{\mathrm{causal}}\right) \overset{%
d}{\longrightarrow }\mathcal{N}\left(0,\Gamma ^{-1}\left(\rho\Delta^{\mathrm{cond}%
}+(1-\rho) \Delta^\mathrm{ehw} \right) \Gamma ^{-1}\right) ,
\end{equation*}

\item[(ii)]
\begin{equation*}
\sqrt{N}\left(\widehat{\theta}_n-\theta_n^{\mathrm{causal},\mathrm{sample}}
\right) \overset{d}{\longrightarrow}\mathcal{N}\left(0,\Gamma ^{-1} \Delta^{%
\mathrm{cond}}\Gamma ^{-1}\right) ,
\end{equation*}

\item[(iii)]
\begin{equation*}
\sqrt{N}\left( \widehat{\theta}_n-\theta_n^{\mathrm{descr}} \right) \overset{%
d}{\longrightarrow }\mathcal{N}\left(0,(1-\rho )\Gamma ^{-1}\Delta^\mathrm{ehw}
\Gamma ^{-1}\right) .
\end{equation*}
\end{itemize}
\end{theorem}

\begin{remark}
For both the population causal and the descriptive estimand the asymptotic
variance in the case with $\rho =0$ reduces to the standard EHW
variance, $\Gamma ^{-1}\Delta^{\mathrm{ehw}}\Gamma ^{-1}$. If the sample size is
non-negligible as a fraction of the population size, $\rho >0$, the
difference between the EHW variance and the finite population
causal variance is positive semi-definite and equal to $\rho \Gamma ^{-1} (\Delta^{%
\mathrm{ehw}}-\Delta^\mathrm{cond})\Gamma ^{-1}$.\hfill$\square$
\end{remark}

\begin{remark}{\sc (The Case with $\rho=0$)}
 The standard setting where we have a random sample from a large population, is covered by the result in Theorem 3, part (i) or  part (iii) with $\rho=0$. For example, when we analyze data from the CPS or PSID, this seems a reasonable perspective. Even if the sampling from the US population is not completely random, it is approximately so, and the sample is certainly small relative to the population. In that case we do not need to worry about  whether we are interested in a causal estimand because the standard methods are valid.
\hfill$\square$
\end{remark}

\begin{remark}{\sc (The Case with $\rho=1$)}
 The case where we observe all units in the population of interest, covered by the result in Theorem 3, part (i) with $\rho=1$, or part (ii)) is also common. For example, we may have all the states in the US, or all the countries in the world, or all the individuals in the population of interest. In that case taking account of the causal nature of the estimand is important, because a descriptive perspective would suggest the standard errors should be zero. 
This  covers the case discussed in  \citet{manski2018right}.
\hfill$\square$
\end{remark}

\begin{remark}{\sc (A Convenience Sample)}
 The setting where we have a 
convenience sample, where the relationship between the sample and the population is murky, is more complicated.  For example, we may have all internet searches during a particular day, or all shopping trips to a single supermarket for a given week. This is in our view an important and common setting.
 In such settings researchers often analyze the data, and report standard errors based on the sampling perspective, as if the sample is a random sample from a large population. Typically they do so implicitly, by simply using standard methods without explicitly describing a sampling process. 
It seems a stretch to view the sample of shopping trips to a particular supermarket on a particular day as a random sample from the population of interest. 
At best it is a systematic sample from the population of interest, e.g., all individuals going to that particular supermarket, rather than a random sample from the population of individuals. However, there is no way to quantify the uncertainty arising from that sampling scheme without data from other periods or other supermarkets.
In that case  we  recommend to analyze the uncertainty relative to the causal sample estimand, and to be clear about what that estimand is in order to provide a conceptually precise measure of uncertainty.
\hfill$\square$
\end{remark}

\begin{remark} {\sc (The case with $\rho\in(0,1)$)}
While  there are clearly many cases in practice where we do observe the entire population of interest, there are also  settings where the following three things hold, at least approximately: $(a)$ the population of interest is finite, $(b)$ the sample is a random sample from this population, and $(c)$, the ratio of sample size to population size $\rho$ is known and large enough for this to matter.  For example,
\citet{muralidharan2017experimentation} discuss a number of randomized experiments in development economics where the study sample was drawn randomly from the population of interest.
\citet{keels2005fifteen} discuss using a 50\% random sample in a mobility study, rather than the full population, for cost or computational reasons.
The Integrated Public Use Microdata Series (IPUMS)  data include a  random sample of 10\% of the census.
Some other recent papers include \citet{dellavigna2017reference}, whose sample consists of a 50\% {\it de facto} random sample of Hungarian citizens older than 14 and younger than 75 in 2002, 
\citet{einav2015response}, who use  a 20 percent random sample of Medicare Part D beneficiaries from 2007 to 2009,
\citet{hanna2014learning}, who randomly selected 117 (from the set of respondents) to participate in an experimental trial, and
\citet{farber2015you}, who uses a random subsample of 2/15 of the drivers in his data set.
Another interesting case, with a more complex sample is
\citet{munnell1996mortgage},   with the sample analyzed containing the population of mortgage applications in the city of Boston in 1990 for Black and Hispanic applicants (1200 obs) and a random sample of White applicants (3300 obs).
\hfill$\square$
\end{remark}

\begin{remark}
Presenting the variance for the general case that includes the case with $\rho=0$ and $\rho=1$  as special cases is helpful because it explicitly connects the two leading perspectives to uncertainty, sampling-based and design-based. It shows that there is no conceptual conflict between our proposed causal perspective and the standard sampling-based perspective on uncertainty, that our perspective merely adds a second source of uncertainty. It also shows that this perspective is particularly relevant when the researcher has observations on the entire population, a case  that previously had not been satisfactorily addressed in the literature.
\hfill$\square$
\end{remark}

\subsection{The Variance Under Correct Specification}

Consider a constant treatment effect assumption, which is required for a
correct specification of a linear regression function as a function that
describes potential outcomes.

\begin{assumption}
\textsc{(Constant Treatment Effects)} \label{assumption:cte}
\begin{equation*}
Y_{n,i}^{*}(u) = u^{\prime }\theta_n+\xi_{n,i},
\end{equation*}
where $\theta_n$ and $\xi_{n,i}$ are non-stochastic.
\end{assumption}

This strengthens Assumption \ref{assumption:linear} by requiring that the $%
\theta_{n,i}$ do not vary by $i$.

Under Assumption \ref{assumption:cte}, Theorem \ref{theorem:causal} implies
that $\theta_n^\mathrm{causal}=\theta_n$ (although it need not be the case
that $\theta^\mathrm{descr}=\theta_n$). Then, for
\begin{equation*}
\lambda_n = \left(\sum_{i=1}^n Z_{n,i}Z_{n,i}^{\prime
}\right)^{-1}\sum_{i=1}^n Z_{n,i}\xi_{n,i}
\end{equation*}
we obtain that equation (\ref{equation:causalres}) holds for $\gamma_n^%
\mathrm{causal}=\Lambda_n^{\prime }\theta_n+\lambda_n$ and $%
\varepsilon_{n,i}=\xi_{n,i}-Z_{n,i}^{\prime }\lambda_n$. In this case, the
residuals, $\varepsilon_{n,i}$, are non-stochastic. As a result, $%
E[X_{n,i}\varepsilon_{n,i}]=E[X_{n,i}]\varepsilon_{n,i}=0$, which implies $%
\Delta^\mu=\Delta^\mathrm{ehw}-\Delta^\mathrm{cond}=0$. This leads to the
following result.

\begin{theorem}
\label{theorem:neyman} Suppose that Assumptions \ref{assumption:assignment}-%
\ref{assumption:cte} hold. Then,
\begin{equation*}
\sqrt{N}\left( \widehat{\theta}_n-\theta^{\mathrm{causal}}_n \right) \overset{d%
}{\longrightarrow }\mathcal{N}\left(0,\Gamma ^{-1}\Delta^\mathrm{ehw}
\Gamma ^{-1}\right) ,
\end{equation*}
irrespective of the value of $\rho$.
\end{theorem}

Notice that the result of the theorem applies also with $\theta_n^{\mathrm{%
causal},\mathrm{sample}}$ replacing $\theta_n^{\mathrm{causal}}$ because the
two parameter vectors are identical (with probability approaching one) under
Assumption \ref{assumption:cte}.

\begin{remark}
The key insight in this theorem is that the asymptotic variance of $%
\widehat\theta_n$ does not depend on the ratio of the sample to the
population size when the regression function is correctly specified.
Therefore, it follows that the usual EHW variance matrix is
correct for $\widehat\theta_n$ under these assumptions. For the special case with $%
X_{n,i}$ binary and no attributes beyond the intercept, this result can be
inferred directly from Neyman's results for randomized experiments (%
\citeauthor{neyman1923}, \citeyear{neyman1923}). In that case, the result of
Theorem \ref{theorem:neyman} follows from the restriction of constant
treatment effects, $Y_{n,i}^{*}(1)-Y_{n,i}^{*}(0)=\theta_n$, which is
extended to the more general case of non-binary regressors in Assumption \ref%
{assumption:cte}. The asymptotic variance of $\widehat{\gamma}_n$, the least
squares estimator of the coefficients on the attributes, still depends on
the ratio of sample to population size, and it can be shown that the
conventional robust EHW estimator continues to over-estimate the
variance of $\widehat{\gamma}_n$. For more details see the earlier version of this paper, \citet{abadieathey}.\hfill$\square$
\end{remark}

\section{Estimating the Variance}

\label{estimatingvariance}

Now let us turn to the problem of estimating the variance for the
descriptive and causal estimands. In what follows, we will use the shorthands $V^{\causal}=\Gamma ^{-1}\left(\rho\Delta^{\mathrm{cond}%
}+(1-\rho) \Delta^\mathrm{ehw} \right) \Gamma ^{-1}$,  $V^{\causal, \sample}=\Gamma ^{-1}\Delta^{\mathrm{cond}} \Gamma ^{-1}$, $V^{\descr}=(1-\rho )\Gamma ^{-1}\Delta^\mathrm{ehw}
\Gamma ^{-1}$, and $V^{\ehw}=\Gamma ^{-1}\Delta^\mathrm{ehw}
\Gamma ^{-1}$. There are four components to the
asymptotic variances, $\rho$, $\Gamma $, $\Delta^\mathrm{ehw}$ and $\Delta^\mathrm{%
cond}$. The first three are straightforward to estimate. $\rho$
can be estimated as $\widehat\rho_n=N/n$, as long as the population size is known. To estimate $\Gamma $, first estimate $\Lambda_n$ as
\begin{equation*}
\widehat{\Lambda}_n=\left(\sum_{i=1}^n R_{n,i} U_{n,i} Z_{n,i}^{\prime
}\right)\left(\sum_{i=1}^n R_{n,i} Z_{n,i} Z_{n,i}^{\prime }\right)^{-1}.
\end{equation*}
Then one can estimate $\Gamma $ as the average of the matrix of outer products
over the sample:
\begin{equation*}
\widehat \Gamma _n=\frac{1}{N}\sum_{i=1}^{n}R_{n,i} \left( U_{n,i}-\widehat{\Lambda%
}_n Z_{n,i} \right) \left( U_{n,i}-\widehat{\Lambda}_n Z_{n,i}
\right)^{\prime }.
\end{equation*}
It is also straightforward to estimate $\Delta^\mathrm{ehw}$. First we
estimate the residuals for the units in the sample, $\widehat{\varepsilon}%
_{n,i}=Y_{n,i}-(U_{n,i}-\widehat{\Lambda}_n Z_{n,i})^{\prime
}\widehat\theta_n -Z_{n,i}^{\prime }\widehat\gamma_n,$ and then we estimate $%
\Delta^\mathrm{ehw}$ as:
\begin{equation*}
\widehat{\Delta}^\mathrm{ehw}_n= \frac{1}{N}\sum_{i=1}^n R_{n,i} (U_{n,i}-%
\widehat{\Lambda}_n Z_{n,i})\,\widehat{\varepsilon}^{\, 2}_{n,i}\, (U_{n,i}-%
\widehat{\Lambda}_n Z_{n,i})^{\prime }.
\end{equation*}
The EHW large sample variance, $V^\mathrm{ehw}$, is then estimated as
\begin{equation*}
\widehat V^\mathrm{ehw}_n=\widehat \Gamma ^{-1}_n\widehat\Delta^\mathrm{ehw}%
_n\widehat \Gamma ^{-1}_n.
\end{equation*}

\begin{lemma}
\label{lemma:convergvar} Suppose Assumptions \ref{assumption:assignment}-\ref%
{assumption:unconf} and \ref{assumption:limits} hold with $\delta=4$. Then,
\begin{equation*}
\widehat V^\mathrm{ehw}_n\overset{p}{\longrightarrow} V^\mathrm{ehw}.
\end{equation*}
\end{lemma}

Let $\widehat V^{\desc}_n=(1-\widehat\rho_n) \widehat V^\mathrm{ehw}_n$. The result of Lemma \ref{lemma:convergvar} immediately implies $\widehat V^{\desc}_n\overset{p}{\rightarrow}
V^{\desc}$.

 It is more challenging to estimate $V^{\causal}$ and $V^{\causal,\sample}$ because they involve $\Delta^\mathrm{cond}$. Estimating $\Delta^{\cond}$ is complicated because of the same reason
that complicates the estimation of the variance
of the average treatment effect estimator in Section %
\ref{section:cases}. In that case there are three terms in the expression
for the variance in equation (\ref{equation:totalvariance}). The first two
are straightforward to estimate, but the third one, $S^2_\theta/n$ cannot be
estimated consistently because we do not observe both potential outcomes for
the same units. Often, researchers use the conservative estimator based on
ignoring $S^2_\theta/n$. If we proceed in the same fashion for the
regression context of Section \ref{main}, we obtain the conservative
estimator $\widehat V^\mathrm{ehw}$, based on ignoring $\Delta^\mu$. We
show, however, that in the presence of attributes we can improve the
variance estimator. We build on \citet{abadie2008estimation}, %
\citet{abadie2014inference}, and \citet{fogarty2016regression} who, in
contexts different than the one studied in this article, have used the
explanatory power of attributes to improve variance estimators.  %
\citet{abadie2008estimation} and \citet{abadie2014inference} do so using 
nearest-neighbor techniques. Here we follow \citet{fogarty2016regression}
and apply linear regression techniques. The proposed estimator replaces the
expectations $E[X_{n,i}\varepsilon_{n,i}]$, which cannot be consistently
estimated, with predictors from a linear least squares projection of
estimates of $X_{n,i}\varepsilon_{n,i}$ on the attributes, $Z_{n,i}$. Let $%
\widehat X_{n,i}= U_{n,i}-\widehat\Lambda_n Z_{n,i}$, and
\begin{equation*}
\widehat G_n = \left(\frac{1}{N}\sum_{i=1}^n R_{n,i} \widehat
X_{n,i}\widehat\varepsilon_{n,i}Z_{n,i}^{\prime }\right)\left(\frac{1}{N}%
\sum_{i=1}^n R_{n,i}Z_{n,i}Z_{n,i}^{\prime }\right)^{-1}.
\end{equation*}
The matrix $\widehat G_n$ contains the coefficients of a least squares
regression of $\widehat X_{n,i}\widehat\varepsilon_{n,i}$ on $Z_{n,i}$. The
next assumption ensures convergence of $\widehat G_n$.

\begin{assumption}
\label{assumption:convergG}
\begin{equation*}
\frac{1}{n}\sum_{i=1}^n E[X_{n,i}\varepsilon_{n,i}]Z_{n,i}^{\prime }
\end{equation*}
has a limit.
\end{assumption}

Consider now the following estimator,
\begin{equation*}
\widehat \Delta^Z_n = \frac{1}{N}\sum_{i=1}^n R_{n,i}\left(\widehat
X_{n,i}\widehat\varepsilon_{n,i}-\widehat G_nZ_{n,i}\right)\left(\widehat
X_{n,i}\widehat\varepsilon_{n,i}-\widehat G_nZ_{n,i}\right)^{\prime }.
\end{equation*}
which uses $\widehat G_nZ_{n,i}$ in lieu of a consistent estimator of $%
E[X_{n,i}\varepsilon_{n,i}]$. Notice that we do not assume that $%
E[X_{n,i}\varepsilon_{n,i}]$ is linear in $Z_{n,i}$. However, we will show
that, as long as the attributes can linearly explain some of the variance in
$\widehat X_{n,i}\widehat\varepsilon_{n,i}$, the estimator $\widehat
\Delta^Z_n$ is smaller (in a matrix sense) than $\widehat\Delta^\mathrm{ehw}%
_n$.
These results are provided in the following lemma.

\begin{lemma}
\label{lemma:convergvar2} Suppose Assumptions \ref{assumption:assignment}-%
\ref{assumption:unconf}, \ref{assumption:limits} and \ref%
{assumption:convergG} hold with $\delta=4$. Then, $0\leq
\widehat\Delta^Z_n\leq \widehat\Delta^\mathrm{ehw}_n$, and $%
\widehat\Delta^Z_n\overset{p}{\rightarrow} \Delta^Z$, where $\Delta^\mathrm{%
cond}\leq \Delta^Z\leq \Delta^\mathrm{ehw}$ (all inequalities are to be
understood in a matrix sense).
\end{lemma}

Estimators of $V^{\causal,\sample}$ and $V^{\causal}$ follow immediately from Lemma \ref{lemma:convergvar2} by
replacing $\Delta^\mathrm{cond}$ with the estimate $\widehat\Delta^Z_n$ in
the asymptotic variance formulas of Theorem \ref{theorem:asympdist}, leading
to $\widehat{V}^{\causal,\sample}_n=\widehat \Gamma ^{-1}_n\widehat\Delta^Z_n\widehat \Gamma ^{-1}_n$ for the estimation of
$V^{\causal, \sample}$ and $\widehat V^{\causal}_n=\widehat\rho_n \widehat{V}^{\causal,\sample}_n + (1-\widehat\rho_n)\widehat V^{\ehw}_n$
for the estimation of $V^{\causal}$.
These estimators are not larger (and typically smaller) than  $\widehat V^{\ehw}_n$ and they remain conservative in
large samples.

\begin{remark}
A special case of the adjusted variance arises when $Z_{n,i}$ is a set of exhaustive and mutually exclusive dummy variables, or if we reduce the information in $Z_{n,i}$ to such indicators. Then, the residuals from regressing $\widehat{X}_{n,i}\widehat{\varepsilon}_{n,i}$ on $Z_{n,i}$ are simply stratum-specific demeaned versions of $\widehat{X}_{n,i}\widehat{\varepsilon}_{n,i}$, and a conservative estimator of $\Delta^{\cond}$ can be obtained using the variance formulas in \cite{wooldridge2001} for standard stratified samples. \hfill $\square$ 
\end{remark}

\section{Simulations}

In this section, we use a simple data-generating process as well as simulations to illustrate the difference between the conventional EHW
variance estimator and the variance estimators proposed in this article. We focus on the case of a
single causal variable, $X_{n,i}$. In addition to the causal variable, the simulations employ an outcome variable, $Y_{n,i}$, and a vector of attributes, $Z_{n,i}$, which consists of a constant equal to one and
$k$ values drawn independently from the standard normal
distribution. The potential outcome function has the form in equation (\ref{equation:pot_outcomes}).
Population values of $\theta_{n,i}$ are generated as independent draws from a normal distribution with mean $Z_{n,i}'\psi$, where $\psi=(0,\psi_1,\ldots,\psi_k)'$ is a $k+1$ vector, and variance $\sigma^2_\theta$.
Population values of $\xi_{n,i}$ and $U_{n,i}$ are generated as independent draws from a normal distribution with mean zero and variance one. Because in this data-generating process $E[U_{n,i}]=0$, it follows that $\Lambda_n$ is a row-vector of zeros and $X_{n,i}=U_{n,i}$.
We use this data-generating process to produce
a population of size $n$. For this data-generating process it can be shown that $\Gamma =1$, $\Delta^{\ehw}=1+3(\psi'\psi+\sigma^2_\theta)$ and $\Delta^{\cond}=1+2(\psi'\psi+\sigma^2_\theta)$, and $\Delta^Z=1+2\psi'\psi+3\sigma^2_{\theta}$ with probability one.
In each simulation repetition,  we sample units at random with probability $\rho$ from the population. As a result, the sample size $N$ is random with $E[N]=n\rho$.
For each sample we estimate $\widehat\theta_n$ by least squares (as in equation (\ref{equation:least_squares})) and a number of variance estimators.

In Table \ref{tabel1} we report the results of the simulations. We consider seven designs. The first column reports the basic design, with $\rho=0.01$ and $n=100000$, so the average sample size is 1000. In this design, there is one stochastic regressor, so $k=1$, and the distribution of the treatment effect, $\theta_{n,i}$, is given by parameter values $\psi=(0,2)'$ and $\sigma^2_\theta=1$. The remaining designs in the second to seventh columns are variations of the basic design in the first column.
In the second design, we increase the dimensionality of $Z_{n,i}$ used for estimation from two to ten. Still, in this design $\psi$ has all entries equal to zero except for $\psi_1=2$, so only the first stochastic regressor matters for the distribution of $\theta_{n,i}$.
In the next design, we change the population size to 10000, so that the average sample size is 100.
In the fourth design, we change the population size to 1000 and the sampling rate, $\rho$,  to one.
In the fifth design, we impose $\psi'\psi=0$, which makes the treatment effect unrelated to the regressors, $Z_{n,i}$. In the sixth design, we set $\sigma^2_\theta=0$, which removes the stochastic part of the treatment effect. In the last design, $\psi'\psi=0$ and $\sigma^2_\theta=0$, so the treatment effect is constant. The first panel of Table \ref{tabel1} provides the parameters of each of the seven simulation designs.

The  second panel of Table \ref{tabel1} reports the standard deviations of $(\widehat\theta_n-\theta_n^\mathrm{descr})$, $(\widehat\theta_n-\theta_n^{\mathrm{causal},\mathrm{sample}})$ and  $(\widehat\theta_n-\theta_n^\mathrm{causal})$ across simulation iterations. The remaining panels report feasible standard errors based on the estimators of Section \ref{estimatingvariance} as well as bootstrap standard errors, along with coverage rates of the corresponding 95 percent confidence intervals. We employ 50000 iterations for the simulations and 1000 bootstrap samples.
The coverage rates in each of the panels of the table are based on the adding and subtracting 1.96 times the standard errors in the first row of that panel.

For the basic design in the first column of Table \ref{tabel1}, $\rho=0.01$ is very small, and EHW and bootstrap standard errors provide accurate estimates of the standard deviations of $(\widehat\theta_n-\theta_n^{\desc})$ and $(\widehat\theta_n-\theta_n^{\causal})$. However, the standard deviation of $(\widehat\theta_n-\theta_n^{\causal,\sample})$ is substantially smaller than that of $(\widehat\theta_n-\theta_n^\desc)$ and
$(\widehat\theta_n-\theta_n^\causal)$, and the EHW and bootstrap variance estimators are very conservative for the sample average causal effect, $\theta^{\causal,\sample}_n$. The variance estimator based on $\widehat V^{\causal,\sample}_n$ is substantially smaller, and still has more than correct coverage for $\theta^{\causal,\sample}$.
Increasing the number of regressors in the second design leaves the result virtually unaffected. The same patterns of results appear in the third column, albeit with less precise variance estimators due to much smaller sample sizes. In the fourth design, $\rho=1$ and, as predicted by the results in section \ref{estimatingvariance}, EHW standard errors greatly overestimate the variability of $(\widehat\theta_n-\theta_n^\mathrm{descr})$. The same is true for bootstrap standard errors. In the fifth, design we go back to the small sampling rate, $\rho=0.01$ and this time $\psi'\psi=0$, so regressors do not explain variation in treatment effects, and $\Delta^Z=\Delta^{\ehw}$. As suggested by the results in sections \ref{main} and \ref{estimatingvariance}, all variance estimators produce similar results in this design.
In the sixth design, where regressors explain all the variation in treatment effects, $\Delta^{Z}=\Delta^{\cond}$ and standard errors based on $\widehat V_n^{\causal,\sample}$ and $\widehat V_n^{\causal}$ closely approximate the standard deviations of $(\widehat\theta_n-\theta_n^{\causal,\sample})$ and $(\widehat\theta_n-\theta_n^{\causal})$, respectively. In the final design with a constant treatment effect, all the variances are similar.

\section{Conclusion}

In this article we study the interpretation of standard errors in regression
analysis when the assumption that the sample is drawn randomly from a much larger
population of interest is not appropriate. We base our results on a potential outcome framework, where the estimands of interest may be
descriptive or causal, and we provide a coherent
interpretation for standard errors that allows for uncertainty coming from
both random sampling and from conditional random assignment. The standard errors estimators proposed in this article
may be different from the conventional ones, and they may vary depending on {\em (i)} the specific nature of the estimand of interest (i.e., descriptive
or causal), {\em (ii)} the fraction of the population represented in the sample, and {\em (iii)} the extent to which measured attributes explain variation in treatment effects.

In the current article we focus exclusively on linear regression models. The concerns we
raise in this article arise in many other settings and for other kinds
of hypotheses, and the implications would need to be worked out for those
settings. Thus, we see this
article as a first step in broader research program.\newpage

\appendix
\centerline{\sc Appendix}

\vskip0.5cm

\centerline{\sc I. A Bayesian Approach}

\vskip0.5cm

Given that we are advocating for a different conceptual approach to modeling
inference, it is useful to look at the problem from more than one
perspective. In this section we consider a Bayesian perspective and
re-analyze the example from Section \ref{section:cases}. Viewing the problem from a Bayesian perspective
reinforces the point that formally modeling the population and the sampling
process leads to the conclusion that inference is different for descriptive
and causal questions. Note that in this discussion the notation will
necessarily be slightly different from the rest of the article; notation and
assumptions introduced in this subsection apply only within this subsection.

Define $\boldsymbol{Y}^*_n(1)$, $\boldsymbol{Y}^*_n(0)$ to be the $n$ vectors
with typical elements $Y^*_i(1)$ and $Y^*_i(0)$, respectively. We view the $n$%
-vectors $\boldsymbol{Y}^*_n(1)$, $\boldsymbol{Y}^*_n(0)$, $\boldsymbol{R}_n$ ,
and $\boldsymbol{X}_n$ as random variables, some observed and some
unobserved. We assume the rows of the $n\times 4$ matrix $[\boldsymbol{Y}^*%
_n(1),\boldsymbol{Y}^*_n(0), \boldsymbol{R}_n, \boldsymbol{X}_n]$ are
exchangeable. Then, by appealing to DeFinetti's theorem, we model this, with
no essential loss of generality (for large $n)$ as the product of $n$
independent and identically distributed random quadruples $%
(Y^*_i(1),Y^*_i(0),R_i,X_i)$ given some unknown parameter $\beta$:
\begin{equation*}
f(\boldsymbol{Y}^*_n(1),\boldsymbol{Y}^*_n(0),\boldsymbol{R}_n,\boldsymbol{X}%
_n)=\prod_{i=1}^{n}f(Y^*_i(1),Y^*_i(0),R_i,X_{i}|\beta).
\end{equation*}
Inference then proceeds by specifying a prior distribution for $\beta $, say
$p(\beta )$. To make this specific, consider the following model. Let $X_i$ and $%
R_i$ have Binomial distributions with parameters $q$ and $\rho$,
\begin{equation*}
\Pr(X_i=1|Y^*_i(1),Y^*_i(0),R_i)=q,\qquad \Pr(R_i=1|Y^*_i(1),Y^*_i(0))=\rho.
\end{equation*}
The pairs $(Y^*_i(1),Y^*_i(0))$ are assumed to be jointly normally distributed:
\begin{equation*}
\left. \left(%
\begin{array}{c}
Y^*_i(1) \\
Y^*_i(0)%
\end{array}
\right)\right|\mu_1 ,\mu_0 ,\sigma^2_1,\sigma^2_0,\kappa\sim \mathcal{N}%
\left(\left(
\begin{array}{c}
\mu_1 \\
\mu_0%
\end{array}
\right), \left(
\begin{array}{cc}
\sigma^2_1 & \kappa\sigma_1\sigma_0 \\
\kappa\sigma_1 \sigma_0 & \sigma^2_0%
\end{array}
\right) \right),
\end{equation*}
so that the full parameter vector is $\beta=(q,\rho,\mu_1,\mu_0
,\sigma^2_1,\sigma^2_0,\kappa)$.

We change the observational scheme slightly from Section \ref{section:cases}
to allow for the analytic derivation of posterior distributions. We assume
that for all units in the population we observe the pair $(R_i,X_i)$, and
for units with $R_i=1$ we observe the outcome $Y_i=Y^*_i(X_i)$. Define $%
\widetilde Y_i=R_i Y_i$, so for all units in the population we observe the
triple $(R_i,X_i,\widetilde Y_i)$. Let $\boldsymbol{R}_n$, $\boldsymbol{X}_n$%
, and $\widetilde{\boldsymbol{Y}}_n$ be the $n$ vectors of these variables. $%
\bar Y_1 $ denotes the average of $Y_i$ in the subpopulation with $R_i=1$
and $X_i=1$, and $\bar Y_0$ denotes the average of $Y_i$ in the
subpopulation with $R_i=1 $ and $X_i=0$.

The descriptive estimand is
\begin{equation*}
\theta_n^\mathrm{descr}=\frac{1}{n_1}\sum_{i=1}^n X_i Y_i-\frac{1}{n_0}%
\sum_{i=1}^n (1-X_i) Y_i.
\end{equation*}
The causal estimand is
\begin{equation*}
\theta_n^\mathrm{causal} =\frac{1}{n} \sum_{i=1}^n \Bigl(Y^*_i(1)-Y^*_i(0)\Bigr).
\end{equation*}
It is interesting to compare these estimands to an additional estimand, the
super-population average treatment effect,
\begin{equation*}
\theta^\mathrm{causal} =\mu_1 -\mu_0.
\end{equation*}
In general these three estimands are distinct, with their own posterior
distributions, but in some cases, notably when $n$ is large, the three
posterior distributions are similar.

It is instructive to consider a very simple case where analytic solutions
for the posterior distribution for $\theta^\mathrm{descr}_n$, $\theta^%
\mathrm{causal}_n$, and $\theta^\mathrm{causal}$ are available. Suppose $%
\sigma^2_1 $, $\sigma^2_0$, $\kappa$ and $q$ are known, so that the only
unknown parameters are the two means $\mu_1$ and $\mu_0$. Finally, let us
use independent, diffuse (improper), prior distributions for $\mu_1$ and $%
\mu_0$.

Then, a standard result is that the posterior distribution for $(\mu_1
,\mu_0)$ given $({{\mathbf{R}}}_n,{\mathbf{X}}_n, \widetilde{{\mathbf{Y}}}%
_n) $ is
\begin{equation*}
\left. \left(
\begin{array}{c}
\mu_1 \\
\mu_0%
\end{array}
\right)\right|{{\mathbf{R}}}_n,{\mathbf{X}}_n,\widetilde{{\mathbf{Y}}}_n
\sim \mathcal{\ N}\left( \left(
\begin{array}{c}
\bar Y_1 \\
\bar Y_0%
\end{array}
\right) , \left(
\begin{array}{cc}
\sigma^2_1/N_1 & 0 \\
0 & \sigma^2_0 /N_0%
\end{array}
\right) \right),
\end{equation*}
where $N_1$ is the number of units with $R_i=1$ and $X_i=1$, and $N_0$ is
the number of units with $R_i=1$ and $X_i=0$. This directly leads to the
posterior distribution for $\theta^\mathrm{causal}$:
\begin{equation*}
\theta^\mathrm{causal} |{{\mathbf{R}}}_n,{\mathbf{X}}_n,\widetilde{{\mathbf{Y%
}}}_n \sim \mathcal{N}\left(\bar Y_1 - \bar Y_0, \frac{\sigma^2_1}{N_1}+%
\frac{\sigma^2_0}{N_0}\right).
\end{equation*}
A longer calculation leads to the posterior distribution for the descriptive
estimand:
\begin{equation*}
\theta^\mathrm{descr}_n |{{\mathbf{R}}}_n,{\mathbf{X}}_n,\widetilde{{\mathbf{%
Y}}}_n \sim \mathcal{N}\left(\bar Y_1 -\bar Y_0, \frac{\sigma^2_1}{N_1}%
\left(1-\frac{N_1}{n_1}\right) + \frac{\sigma^2_0}{N_0}\left(1-\frac{N_0}{n_0%
}\right) \right).
\end{equation*}
The implied posterior interval for $\theta^\mathrm{descr}_n $ is very
similar to the corresponding confidence interval based on the normal
approximation to the sampling distribution for $\bar Y_1 -\bar Y_0$. If $n_1$
and $n_0$ are large, this posterior distribution is close to the posterior
distribution of the causal estimand. If, on the other hand, $N_1=n_1$ and $%
N_0=n_0$, then the posterior distribution of the descriptive estimand
becomes degenerate and centered at $\bar Y_1 -\bar Y_0$.

A somewhat longer calculation for $\theta^\mathrm{causal}_n$ leads to
\begin{equation*}
\theta^\mathrm{causal}_n |{{\mathbf{R}}}_n,{\mathbf{X}}_n,\widetilde{{%
\mathbf{Y}}}_n \sim \mathcal{N}\left(\bar Y_1 -\bar Y_0 ,\frac{N_0}{n^2}%
\sigma^2_1 (1-\kappa^2) +\frac{N_1}{n^2}\sigma^2_0(1-\kappa^2) \right.
\end{equation*}
\begin{equation*}
\hskip2cm +\frac{n-N}{n^2}\sigma^2_1 +\frac{n-N}{n^2}\sigma^2_0 -2\frac{n-N}{%
n^2}\kappa\sigma_1 \sigma_0
\end{equation*}
\begin{equation*}
\hskip2cm \left. +\frac{\sigma^2_1}{N_1} \left(1-\left(1-\kappa\frac{\sigma_0%
}{\sigma_1}\right)\frac{N_1}{n}\right)^2 +\frac{\sigma^2_0}{N_0}%
\left(1-\left(1-\kappa\frac{\sigma_1}{\sigma_0}\right)\frac{N_0}{n}\right)^2
\right).
\end{equation*}
Consider the special case of constant treatment effects, where $%
Y_i(1)-Y_i(0)=\mu_1-\mu_0$. Then, $\kappa=1$, and $\sigma_1=\sigma_0$, and
the posterior distribution of $\theta^\mathrm{causal}_n$ is the same as the
posterior distribution of $\theta^\mathrm{causal}$. The same posterior
distribution arises in the limit if $n$ goes to infinity, regardless of the
values of $\kappa$, $\sigma_1$, and $\sigma_0$.

To summarize, if the population is large, relative to the sample, the posterior
distributions of $\theta^\mathrm{descr}_n$, $\theta^\mathrm{causal}_n$ and $%
\theta^\mathrm{causal}$ agree. However, if the population is small, the
three posterior distributions differ, and the researcher needs to be precise
in defining the estimand. In such cases, simply focusing on the
super-population estimand $\theta^\mathrm{causal} =\mu_1 -\mu_0$ is arguably
not appropriate, and the posterior inferences for such estimands will differ
from those for other estimands such as $\theta^\mathrm{causal}_n$ or $\theta^%
\mathrm{descr}_n$. \vskip0.5cm

\centerline{\sc II. Proofs}

\bigskip

\textbf{Proof of Lemma \ref{lemma:convergence2moments}}: See supplementary
appendix.\hfill$\square$ \bigskip

\textbf{Proof of Theorem \ref{theorem:causal}}: For $n$ large enough $%
\sum_{i=1}^n Z_{n,i}Z_{n,i}^{\prime }$ is full rank and $\Lambda_n$ exists,
so $\Omega^{ZX}_n=0$. This implies
\begin{equation*}
\theta^{\mathrm{causal}}_n = \left(\sum_{i=1}^n E[X_{n,i}X_{n,i}^{\prime
}]\right)^{-1}\sum_{i=1}^n E[X_{n,i}Y_{n,i}].
\end{equation*}
Moreover, for $n$ large enough, $\Lambda_n=B_n$, which implies $E[X_{n,i}]=0$%
, $\widetilde \Omega^{XZ}_n=0$, and
\begin{equation*}
\theta^{\mathrm{causal},\mathrm{sample}}_n = \left(\sum_{i=1}^n
R_{n,i}E[X_{n,i}X_{n,i}^{\prime }]\right)^{-1}\sum_{i=1}^n
R_{n,i}E[X_{n,i}Y_{n,i}]
\end{equation*}
with probability approaching one. Now,
\begin{align*}
E[X_{n,i}Y_{n,i}]&=E[X_{n,i}U_{n,i}^{\prime
}]\theta_{n,i}+E[X_{n,i}]\xi_{n,i} \\
&=E[X_{n,i}X_{n,i}^{\prime }]\theta_{n,i}.
\end{align*}
implies the results.\hfill$\square$ \bigskip

{\bf Proof of Theorem \ref{theorem:nonlinear}}: Let $\nabla Y^*_{n,i}()$ be the gradient of $Y^*_{n,i}()$. By the mean value theorem there exist sets $\mathcal T_{n,i}\subseteq [0,1]$ such that for any $t_{n,i}\in\mathcal T_{n,i}$, we have $Y^*_{n,i}(U_{n,i})=Y^*_{n,i}(B_nZ_{n,i})+X_{n,i}'\nabla Y^*_{n,i}(B_nZ_{n,i}+t_{n,i}X_{n,i})$. We define $\varphi_{n,i}=\nabla Y^*_{n,i}(v_{n,i})$, where $v_{n,i}=B_nZ_{n,i}+\bar t_{n,i}X_{n,i}$ and $\bar t_{n,i}=\sup\,\mathcal T_{n,i}$. Now, $E[X_{n,i}Y_{n,i}]=E[X_{n,i}]Y^*_{n,i}(B_nZ_{n,i})+E[X_{n,i}'\varphi_{n,i}]=E[X_{n,i}'\varphi_{n,i}]$. The rest of the proof is as for Theorem \ref{theorem:causal}.\hfill
$\square$
\bigskip

The following lemma will be useful for establishing asymptotic normality.

\begin{lemma}
\label{lemma:twaalf} Let $V_{n,i}$ is a row-wise independent triangular
array and $\mu_{n,i}=E[V_{n,i}]$. Suppose that $R_{n,1},\ldots, R_{n,n}$ are
independent of $V_{n,1},\ldots, V_{n,n}$ and that Assumption \ref%
{assumption:randomsampling} holds. Moreover, assume that
\begin{equation*}
\frac{1}{n}\sum_{i=1}^n E\left[|V_{n,i}|^{2+\delta}\right]
\end{equation*}
is bounded for some $\delta>0$,
\begin{equation}  \label{equation:zerosum}
\sum_{i=1}^n \mu_{n,i}= 0,
\end{equation}
\begin{equation*}
\frac{1}{n}\sum_{i=1}^n \mbox{\em var}(V_{n,i})\rightarrow \sigma^2,
\end{equation*}
and
\begin{equation*}
\frac{1}{n}\sum_{i=1}^n \mu_{n,i}^2\rightarrow \kappa^2,
\end{equation*}
where $\sigma^2+(1-\rho)\kappa^2>0$. Then
\begin{equation*}
\frac{1}{\sqrt N}\sum_{i=1}^n R_{n,i} V_{n,i}\overset{d}{\longrightarrow}
\mathcal{N }(0,\sigma^2+(1-\rho)\kappa^2),
\end{equation*}
where $N=\sum_{i=1}^n R_{n,i}$.
\end{lemma}

\textbf{Proof:} Notice that
\begin{equation*}
E\left[\frac{N}{n\rho_n}\right]=1
\end{equation*}
and
\begin{equation*}
\mbox{var}\left(\frac{N}{n\rho_n}\right)=\frac{n\rho_n(1-\rho_n)}{(n\rho_n)^2%
}\rightarrow 0.
\end{equation*}
Now the continuous mapping theorem implies
\begin{equation*}
\left(\frac{n\rho_n}{N}\right)^{1/2}\overset{p}{\longrightarrow} 1.
\end{equation*}
As a result, it is enough to prove
\begin{equation*}
\frac{1}{\sqrt n}\sum_{i=1}^n \frac{R_{n,i}}{\sqrt{\rho_n}}
V_{n,i}\rightarrow \mathcal{N }(0,\sigma^2+(1-\rho)\kappa^2).
\end{equation*}
Let
\begin{equation*}
s_n^2 = \frac{1}{n}\sum_{i=1}^n \left(\mbox{var}(V_{n,i})+(1-\rho_n)%
\mu_{n,i}^2\right).
\end{equation*}
Consider $n$ large enough so $s_n^2>0$. Notice that, for $i=1,\ldots, n$,
\begin{equation*}
E\left[\frac{R_{n,i}V_{n,i}-\rho_n\mu_{n,i}}{s_n\sqrt{n\rho_n}}\right]=0,
\end{equation*}
and
\begin{align*}
\mbox{var}\left(R_{n,i}V_{n,i}-\rho_n\mu_{n,i}\right)&=\rho_nE[V_{n,i}^2]-%
\rho_n^2\mu_{n,i}^2 \\
&=\rho_n\left(\mbox{var}(V_{n,i})+(1-\rho_n)\mu_{n,i}^2\right).
\end{align*}
Therefore,
\begin{align*}
\sum_{i=1}^n \mbox{var}\left(\frac{R_{n,i}V_{n,i}-\rho_n\mu_{n,i}}{s_n\sqrt{%
n\rho_n}}\right)=1.
\end{align*}
Using $\rho_n\leq \rho_n^{1/(2+\delta)}$, $|\mu_{n,i}|^{2+\delta}\leq
E[|V_{n,i}|^{2+\delta}]$, and Minkowski's inequality, we obtain:
\begin{align*}
\sum_{i=1}^n E\left[\left|\frac{R_{n,i}V_{n,i}-\rho_n\mu_{n,i}}{s_n\sqrt{%
n\rho_n}}\right|^{2+\delta}\right] &\leq \frac{1}{s_n^{2+\delta}(n%
\rho_n)^{1+\delta/2}}\sum_{i=1}^n\left(\rho_n^{\frac{1}{2+\delta}}\left(E%
\left[|V_{n,i}|^{2+\delta}\right]\right)^{\frac{1}{2+\delta}%
}+\rho_n|\mu_{n,i}|\right)^{2+\delta} \\
&\leq \frac{2^{2+\delta}\rho_n}{s_n^{2+\delta}(n\rho_n)^{1+\delta/2}}%
\sum_{i=1}^n E\left[|V_{n,i}|^{2+\delta}\right] \\
&=\frac{2^{2+\delta}}{s_n^{2+\delta}(n\rho_n)^{\delta/2}}\left(\frac{1}{n}%
\sum_{i=1}^n E\left[|V_{n,i}|^{2+\delta}\right]\right)\rightarrow 0.
\end{align*}
Applying Liapunov's theorem (see, e.g., \citeauthor{davidson1994stochastic}, %
\citeyear{davidson1994stochastic}), we obtain
\begin{equation*}
\sum_{i=1}^n \frac{R_{n,i}V_{n,i}-\rho_n\mu_{n,i}}{s_n\sqrt{n\rho_n}}\overset%
{d}{\longrightarrow} \mathcal{N }(0,1).
\end{equation*}
Now, the result of the lemma follows from equation (\ref{equation:zerosum})
and from $s_n/\sqrt{\sigma^2+(1-\rho)\kappa^2}\rightarrow 1$.\hfill$\square$
\bigskip

\begin{lemma}
\label{lemma:acht} Suppose Assumptions \ref{assumption:assignment}-\ref%
{assumption:limits} hold, and let $\Delta^\mu = \Delta^\mathrm{ehw}-\Delta^%
\mathrm{cond}$, $\widetilde\varepsilon_{n,i}=Y_{n,i}-X_{n,i}^{\prime
}\theta_n^{\mathrm{causal},\mathrm{sample}}-X_{n,i}^{\prime }\gamma_n^{%
\mathrm{causal},\mathrm{sample}}$, and $\nu_{n,i}=Y_{n,i}-X_{n,i}^{\prime
}\theta_n^{\mathrm{descr}}-X_{n,i}^{\prime }\gamma_n^{\mathrm{descr}}$. Then,

\begin{itemize}
\setlength\itemsep{-2pt}
\item[(i)]
\begin{equation*}
\frac{1}{\sqrt{N}}\sum_{i=1}^n R_{n,i} X_{n,i} \varepsilon_{n,i} \overset{d}{%
\longrightarrow }{\mathcal{N}}(0, \Delta^{\mathrm{cond}}+(1-\rho)
\Delta^\mu),
\end{equation*}

\item[(ii)]
\begin{equation*}
\frac{1}{\sqrt{N}}\sum_{i=1}^n R_{n,i} X_{n,i} \widetilde\varepsilon _{n,i}
\overset{d}{\longrightarrow }{\mathcal{N}}(0,\Delta^{\mathrm{cond}}),
\end{equation*}

\item[(iii)]
\begin{equation*}
\frac{1}{\sqrt{N}}\sum_{i=1}^n R_{n,i} X_{n,i}\nu_{n,i} \overset{d}{%
\longrightarrow }{\mathcal{N}}(0,(1-\rho)\Delta ^{\mathrm{ehw}}).
\end{equation*}
\end{itemize}
\end{lemma}

\noindent\textbf{Proof of Lemma \ref{lemma:acht}:} To prove $(i)$, consider $%
V_{n,i}=a^{\prime }X_{n,i}\varepsilon_{n,i}$ for $a\in \mathbb{R}^k$. We
will verify the conditions Lemma \ref{lemma:twaalf}. Notice that,
\begin{align*}
\frac{1}{n}\sum_{i=1}^n E\left[|V_{n,i}|^{2+\delta}\right]&\leq \frac{%
\|a\|^{2+\delta}}{n}\sum_{i=1}^n E\left[\|X_{n,i}\|^{2+\delta}\big(%
|Y_{n,i}|+ \|X_{n,i}\|\|\theta_n\|+\|Z_{n,i}\|\|\gamma_n\|\big)^{2+\delta}%
\right].
\end{align*}
By Minkowski's inequality and Assumption \ref{assumption:moments}, the
right-hand side of last equation is bounded. In addition,
\begin{equation*}
\sum_{i=1}^n \mu_{n,i}=a^{\prime }\sum_{i=1}^n E[X_{n,i}\varepsilon_{n,i}]=0.
\end{equation*}
Let $a\neq 0$. Then,
\begin{align*}
\frac{1}{n}\sum_{i=1}^n \mbox{var}(V_{n,i})=a^{\prime }\left(\frac{1}{n}%
\sum_{i=1}^n\mbox{var}\left(X_{n,i}\varepsilon_{n,i}\right)\right)a%
\rightarrow a^{\prime cond} a>0.
\end{align*}
\begin{equation*}
\frac{1}{n}\sum_{i=1}^n \mu_{n,i}^2 = a^{\prime }\left(\frac{1}{n}%
\sum_{i=1}^n E[X_{n,i}\varepsilon_{n,i}]E[\varepsilon_{n,i}X_{n,i}^{\prime
}]\right)a\rightarrow a^{\prime \mu} a.
\end{equation*}
This implies
\begin{equation*}
a^{\prime }\left(\frac{1}{\sqrt N}\sum_{i=1}^n
R_{n,i}X_{n,i}\varepsilon_{n,i}\right)\overset{d}{\rightarrow} \mathcal{N }%
(0,a^{\prime cond} + (1-\rho) \Delta^{\mu}) a).
\end{equation*}
Using the Cramer-Wold device, this implies
\begin{equation*}
\frac{1}{\sqrt N}\sum_{i=1}^n R_{n,i}X_{n,i}\varepsilon_{n,i}\overset{d}{%
\rightarrow} \mathcal{N }(0,\Delta^{cond} + (1-\rho) \Delta^{\mu}).
\end{equation*}
The proofs of $(ii)$ and $(iii)$ are similar.\hfill$\square$\bigskip

\textbf{Proof of Theorem \ref{theorem:asympdist} }: To prove $(i)$, notice
that
\begin{equation*}
\sum_{i=1}^{n}R_{n,i}\left(
\begin{array}{cc}
X_{n,i}X_{n,i}^{\prime } & X_{n,i}Z_{n,i}^{\prime } \\
Z_{n,i}X_{n,i}^{\prime } & Z_{n,i}Z_{n,i}^{\prime }%
\end{array}%
\right)
\end{equation*}%
is invertible with probability approaching one. Then,
\begin{align*}
\left(
\begin{array}{c}
\widehat{\theta }_{n} \\
\widehat{\gamma }_{n}%
\end{array}%
\right) & =\left( \sum_{i=1}^{n}R_{n,i}\left(
\begin{array}{cc}
X_{n,i}X_{n,i}^{\prime } & X_{n,i}Z_{n,i}^{\prime } \\
Z_{n,i}X_{n,i}^{\prime } & Z_{n,i}Z_{n,i}^{\prime }%
\end{array}%
\right) \right) ^{-1}\sum_{i=1}^{n}R_{n,i}\left(
\begin{array}{c}
X_{n,i}Y_{n,i} \\
Z_{n,i}Y_{n,i}%
\end{array}%
\right) \\
& =\left(
\begin{array}{c}
\theta _{n}^{causal} \\
\gamma _{n}^{causal}%
\end{array}%
\right) +\left( \sum_{i=1}^{n}R_{n,i}\left(
\begin{array}{cc}
X_{n,i}X_{n,i}^{\prime } & X_{n,i}Z_{n,i}^{\prime } \\
Z_{n,i}X_{n,i}^{\prime } & Z_{n,i}Z_{n,i}^{\prime }%
\end{array}%
\right) \right) ^{-1}\sum_{i=1}^{n}R_{n,i}\left(
\begin{array}{c}
X_{n,i}\varepsilon _{n,i} \\
Z_{n,i}\varepsilon _{n,i}%
\end{array}%
\right) .
\end{align*}%
Therefore,
\begin{align*}
\sqrt{N}\left(
\begin{array}{c}
\widehat{\theta }_{n}-\theta _{n}^{causal} \\
\widehat{\gamma }_{n}-\gamma _{n}^{causal}%
\end{array}%
\right) & =\left( \frac{1}{N}\sum_{i=1}^{n}R_{n,i}\left(
\begin{array}{cc}
X_{n,i}X_{n,i}^{\prime } & X_{n,i}Z_{n,i}^{\prime } \\
Z_{n,i}X_{n,i}^{\prime } & Z_{n,i}Z_{n,i}^{\prime }%
\end{array}%
\right) \right) ^{-1}\frac{1}{\sqrt{N}}\sum_{i=1}^{n}R_{n,i}\left(
\begin{array}{c}
X_{n,i}\varepsilon _{n,i} \\
Z_{n,i}\varepsilon _{n,i}%
\end{array}%
\right) \\
& =\left(
\begin{array}{cc}
\Omega _{n}^{XX} & \Omega _{n}^{XZ} \\
\Omega _{n}^{ZX} & \Omega _{n}^{ZZ}%
\end{array}%
\right) ^{-1}\frac{1}{\sqrt{N}}\sum_{i=1}^{n}R_{n,i}\left(
\begin{array}{c}
X_{n,i}\varepsilon _{n,i} \\
Z_{n,i}\varepsilon _{n,i}%
\end{array}%
\right) +r_{n},
\end{align*}%
where
\begin{equation*}
r_{n}=\left[ \left(
\begin{array}{cc}
\widetilde{W}_{n}^{XX} & \widetilde{W}_{n}^{XZ} \\
\widetilde{W}_{n}^{ZX} & \widetilde{W}_{n}^{ZZ}%
\end{array}%
\right) ^{-1}-\left(
\begin{array}{cc}
\Omega _{n}^{XX} & \Omega _{n}^{XZ} \\
\Omega _{n}^{ZX} & \Omega _{n}^{ZZ}%
\end{array}%
\right) ^{-1}\right] \frac{1}{\sqrt{N}}\sum_{i=1}^{n}R_{n,i}\left(
\begin{array}{c}
X_{n,i}\varepsilon _{n,i} \\
Z_{n,i}\varepsilon _{n,i}%
\end{array}%
\right) .
\end{equation*}%
Because (i) $\Omega _{n}^{XZ}=0$, (ii) the first term of $r_{n}$ is $%
o_{p}(1) $, and (iii) $(1/\sqrt{N})\sum_{i=1}^{n}R_{n,i}X_{n,i}\varepsilon
_{n,i}$ is $O_{p}(1)$ (under the conditions stated above), it follows that%
\begin{equation*}
\sqrt{N}(\widehat{\theta }_{n}-\theta _{n}^{causal})=\left( \Omega
_{n}^{XX}\right) ^{-1}\frac{1}{\sqrt{N}}\sum_{i=1}^{n}R_{n,i}X_{n,i}%
\varepsilon _{n,i}+o_{p}(1)
\end{equation*}%
if we can show%
\begin{equation*}
\left( 1/\sqrt{N}\right) \sum_{i=1}^{n}R_{n,i}Z_{n,i}\varepsilon
_{n,i}=O_{p}(1).
\end{equation*}%
We can write this standardized sum as

\begin{equation*}
\left( n\rho _{n}/N\right) ^{1/2}\left[ n^{-1/2}\sum_{i=1}^{n}\left( R_{n,i}/%
\sqrt{\rho _{n}}\right) Z_{n,i}\varepsilon _{n,i}\right] .
\end{equation*}%
As shown in Lemma A.1, $\left( n\rho _{n}/N\right) ^{1/2}\overset{p}{%
\rightarrow }1$. Therefore, it suffices to show

\begin{equation*}
n^{-1/2}\sum_{i=1}^{n}\left( R_{n,i}/\sqrt{\rho _{n}}\right)
Z_{n,i}\varepsilon _{n,i}=O_{p}(1).
\end{equation*}%
This expression has zero mean because $R_{n,i}$ is independent of $%
\varepsilon _{n,i}$ and

\begin{equation*}
\sum_{i=1}^{n}Z_{n,i}\mathrm{E}\left( \varepsilon _{n,i}\right) =0.
\end{equation*}%
We can study each element of the vector separately. By Chebyshev's
inequality it suffices to show that the variances are bounded. Consider the $%
j^{th}$ element. Then, by independence across $i$,

\begin{eqnarray*}
\mathrm{var}\left[ n^{-1/2}\sum_{i=1}^{n}\left( R_{n,i}/\sqrt{\rho _{n}}%
\right) Z_{n,i,j}\varepsilon _{n,i}\right] &=&n^{-1}\sum_{i=1}^{n}\mathrm{var%
}\left[ \left( R_{n,i}/\sqrt{\rho _{n}}\right) Z_{n,i,j}\varepsilon _{n,i}%
\right] \\
&\leq &n^{-1}\sum_{i=1}^{n}E\left\{ \left[ \left( R_{n,i}/\sqrt{%
\rho _{n}}\right) Z_{n,i,j}\varepsilon _{n,i}\right] \right\} ^{2} \\
&=&n^{-1}\sum_{i=1}^{n}Z_{n,i,j}^{2}E\left( \varepsilon
_{n,i,j}^{2}\right)
\end{eqnarray*}%
where the last equality holds because $E[R_{n,i}] =\rho
_{n}$ and $Z_{n,i,j}$ is nonrandom. Both $Z_{n,i,j}^{2}$ and $E[\varepsilon _{n,i}^{2}]$ are bounded by Assumption
\ref{assumption:moments}, and so this completes the proof. The proofs of $(ii)$
and $(iii)$ are analogous.\newline
\strut \hfill $\square $

\noindent\textbf{Proof of Theorem \ref{theorem:neyman}:} The result follows
directly $E[X_{n,i}\varepsilon_{n,i}]=0$. \hfill$\square$

\noindent\textbf{Proof of Lemma \ref{lemma:convergvar}:} First, notice that
(with probability approaching one) $\Lambda_n$ exists and it is equal to $%
B_n $. This implies,
\begin{align*}
\widehat\Lambda_n -\Lambda_n= \left(\frac{1}{N}\sum_{i=1}^n
R_{n,i}X_{n,i}Z_{n,i}^{\prime }\right)\left(\frac{1}{N}\sum_{i=1}^n
R_{n,i}Z_{n,i}Z_{n,i}^{\prime }\right)^{-1}
\end{align*}
which converges to zero in probability by Lemma \ref%
{lemma:convergence2moments} and Assumption \ref{assumption:convergence}.
Direct calculations yield
\begin{align*}
\widehat \Gamma _n -\widetilde W^{XX}_n &= (\widehat\Lambda_n
-\Lambda_n)\widetilde W^{ZZ}_n(\widehat\Lambda_n -\Lambda_n)^{\prime
}-\widetilde W^{XZ}_n(\widehat\Lambda_n -\Lambda_n)^{\prime
}-(\widehat\Lambda_n -\Lambda_n)\widetilde W^{XZ}_n\overset{p}{\rightarrow}
0.
\end{align*}
Now, Lemma \ref{lemma:convergence2moments} and Assumption \ref%
{assumption:convergence} imply $\widehat \Gamma _n\overset{p}{\rightarrow} \Gamma $,
where $\Gamma $ is full rank. Theorem \ref{theorem:asympdist} diretclty implies $%
\widehat\theta_n-\theta^\mathrm{causal}_n\overset{p}{\rightarrow} 0$. $%
\widehat\gamma_n-\gamma^\mathrm{causal}_n\overset{p}{\rightarrow} 0$ follows
from Lemma \ref{lemma:convergence2moments}. Let
\begin{align*}
\breve\Delta^\mathrm{ehw}_n&=\frac{1}{N}\sum_{i=1}^n
R_{n,i}X_{n,i}\widehat\varepsilon_{n,i}^2X_{n,i}^{\prime },\qquad
\widetilde\Delta^\mathrm{ehw}_n=\frac{1}{N}\sum_{i=1}^n
R_{n,i}X_{n,i}\varepsilon_{n,i}^2X_{n,i}^{\prime }, \\
\shortintertext{and} \Delta^\mathrm{ehw}_n&=\frac{1}{n}\sum_{i=1}^n
E[X_{n,i}\varepsilon_{n,i}^2X_{n,i}^{\prime }].
\end{align*}
Let $\alpha$ be a multi-index of dimension equal to the length of $%
T_{n,i}=(Y_{n,i} : X_{n,i}^{\prime }: Z_{n,i}^{\prime })$. In addition, let
\begin{align*}
\widetilde T_n^\alpha&=\frac{1}{N}\sum_{i=1}^n \widetilde T_{n,i}^\alpha =
\frac{1}{N}\sum_{i=1}^n R_{n,i}T_{n,i}^\alpha, \\
\shortintertext{and} \Psi_n^\alpha&=\frac{1}{n}\sum_{i=1}^n
E[W_{n,i}^\alpha].
\end{align*}
Using the same argument as in the proof of Lemma \ref%
{lemma:convergence2moments} and given that Assumption \ref%
{assumption:moments} holds with $\delta=4$, it follows that $\widetilde
T_n^\alpha-\Psi_n^\alpha\overset{p}{\rightarrow} 0$ for $|\alpha|\leq 4$.
This result directly implies $\widetilde\Delta^\mathrm{ehw}_n-\Delta^\mathrm{%
ehw}_n\overset{p}{\rightarrow} 0$. By the same argument plus convergence of $%
\widehat\theta_n$ and $\widehat\gamma_n$, it follows that $\widehat\Delta^%
\mathrm{ehw}_n-\breve\Delta^\mathrm{ehw}_n\overset{p}{\rightarrow} 0$ and $%
\breve\Delta^\mathrm{ehw}_n-\widetilde\Delta^\mathrm{ehw}_n\overset{p}{%
\rightarrow} 0$. Now, the result follows from $\widehat\Delta^\mathrm{ehw}%
_n-\Delta^\mathrm{ehw}=(\widehat\Delta^\mathrm{ehw}_n-\breve\Delta^\mathrm{%
ehw}_n)+(\breve\Delta^\mathrm{ehw}_n-\widetilde\Delta^\mathrm{ehw}%
_n)+(\widetilde\Delta^\mathrm{ehw}_n-\Delta^\mathrm{ehw}_n)+(\Delta^\mathrm{%
ehw}_n-\Delta^\mathrm{ehw})\overset{p}{\rightarrow} 0$, where the last
difference goes to zero by Assumption \ref{assumption:limits}.\hfill$\square$

\noindent\textbf{Proof of Lemma \ref{lemma:convergvar2}:} Notice that,
\begin{align*}
\widehat \Delta^Z_n &= \widehat \Delta^\mathrm{ehw}_n-\widehat\Delta_n^%
\mathrm{proj}, \ \ \ \mathrm{where} \ \ \widehat\Delta_n^\mathrm{proj} =
\frac{1}{N}\sum_{i=1}^n R_{n,i}\widehat G_n Z_{n,i}Z_{n,i}^{\prime }\widehat
G_n^{\prime },
\end{align*}
so that $\widehat \Delta^Z_n$ is no larger than $\widehat \Delta^\mathrm{ehw}%
_n$ in a matrix sense.

Let
\begin{equation*}
G_n = \left(\frac{1}{n}\sum_{i=1}^n
E[X_{n,i}\varepsilon_{n,i}]Z_{n,i}^{\prime }\right)\left(\frac{1}{n}%
\sum_{i=1}^n Z_{n,i}Z_{n,i}^{\prime }\right)^{-1},
\end{equation*}
be the expected value of $\widehat G_n$. Under the assumptions of Lemma \ref%
{lemma:convergvar} and using the same argument as in the proof of that
lemma, we obtain $\widehat G_n-G_n\overset{p}{\rightarrow} 0$. Therefore, $%
\widehat \Delta^\mathrm{proj}_n -\Delta^\mathrm{proj}_n\overset{p}{%
\rightarrow} 0$, where
\begin{equation*}
\Delta_n^\mathrm{proj} = \frac{1}{n}\sum_{i=1}^n G_n Z_{n,i}Z_{n,i}^{\prime
}G_n^{\prime }.
\end{equation*}
Moreover, $\widehat \Delta^Z_n -\Delta^Z_n\overset{p}{\rightarrow} 0$, where
$\Delta^Z_n = \Delta^\mathrm{ehw}_n-\Delta^\mathrm{proj}_n$ and
\begin{equation*}
\Delta_n^\mathrm{ehw} = \frac{1}{n}\sum_{i=1}^n
E[X_{n,i}\varepsilon_{n,i}^2X_{n,i}^{\prime }].
\end{equation*}
Let
\begin{equation*}
\Delta_n^\mu = \frac{1}{n}\sum_{i=1}^n
E[X_{n,i}\varepsilon_{n,i}]E[\varepsilon_{n,i}X_{n,i}^{\prime }].
\end{equation*}
Notice that
\begin{align*}
\Delta_n^\mu-\Delta_n^\mathrm{proj}&= \frac{1}{n}\sum_{i=1}^n
E[X_{n,i}\varepsilon_{n,i}]E[\varepsilon_{n,i}X_{n,i}^{\prime }] \\
&-\left(\frac{1}{n}\sum_{i=1}^n E[X_{n,i}\varepsilon_{n,i}]Z_{n,i}^{\prime
}\right)\left(\frac{1}{n}\sum_{i=1}^n Z_{n,i}Z_{n,i}^{\prime
}\right)^{-1}\left(\frac{1}{n}\sum_{i=1}^n Z_{n,i}E[\varepsilon_{n,i}
X_{n,i}^{\prime }]\right).
\end{align*}
Let $\boldsymbol{A}_n$ and $\boldsymbol{D}_n$ be the matrices with $i$-th
rows equal to $E[\varepsilon_{n,i}X_{n,i}^{\prime }]/\sqrt n$ and $%
Z_{n,i}^{\prime }/\sqrt n$, respectively. Let $\boldsymbol{I}_n$ be the
identity matrix of size $n$. Then,
\begin{equation*}
\Delta_n^\mu-\Delta_n^\mathrm{proj}=\boldsymbol{A}_n^{\prime }(\boldsymbol{I}%
_n-\boldsymbol{D}_n(\boldsymbol{D}_n^{\prime }\boldsymbol{D}_n)^{-1}%
\boldsymbol{D}_n^{\prime })\boldsymbol{A}_n,
\end{equation*}
which is positive semi-definite. Because $\Delta^\mathrm{cond}_n=\Delta^%
\mathrm{ehw}_n-\Delta^\mu_n$, we obtain,
\begin{equation*}
\Delta^\mathrm{cond}_n\leq \Delta^Z_n\leq \Delta^\mathrm{ehw}_n
\end{equation*}
where the inequalities are to be understood in a matrix sense. Now, it
follow from Assumption \ref{assumption:convergG} that $G_n$ and, therefore, $%
\Delta^\mathrm{proj}_n$ and $\Delta^Z_n$ have limits. Then,
\begin{equation*}
\Delta^\mathrm{cond}\leq \Delta^Z\leq \Delta^\mathrm{ehw}
\end{equation*}
where $\Delta^\mathrm{cond}$, $\Delta^Z$, and $\Delta^\mathrm{ehw}$ are the
limits of $\Delta^\mathrm{cond}_n$, $\Delta^Z_n$, and $\Delta^\mathrm{ehw}_n$%
, respectively.\hfill$\square$

\bibliographystyle{chicago}
\bibliography{references}

\begin{thebibliography}{}

\bibitem[\protect\citeauthoryear{Abadie, Athey, Imbens, and Wooldridge}{Abadie
  et~al.}{2017}]{abadie2017should}
Abadie, A., S.~Athey, G.~W. Imbens, and J.~Wooldridge (2017).
\newblock When should you adjust standard errors for clustering?
\newblock Technical report, National Bureau of Economic Research.

\bibitem[\protect\citeauthoryear{Abadie, Athey, Imbens, and Wooldridge}{Abadie
  et~al.}{2014}]{abadieathey}
Abadie, A., S.~Athey, G.~W. Imbens, and J.~M. Wooldridge (2014).
\newblock Finite population causal standard errors.
\newblock Technical report, National Bureau of Economic Research.

\bibitem[\protect\citeauthoryear{Abadie and Imbens}{Abadie and
  Imbens}{2008}]{abadie2008estimation}
Abadie, A. and G.~W. Imbens (2008).
\newblock Estimation of the conditional variance in paired experiments.
\newblock {\em Annales d'Economie et de Statistique\/}, 175--187.

\bibitem[\protect\citeauthoryear{Abadie, Imbens, and Zheng}{Abadie
  et~al.}{2014}]{abadie2014inference}
Abadie, A., G.~W. Imbens, and F.~Zheng (2014).
\newblock Inference for misspecified models with fixed regressors.
\newblock {\em Journal of the American Statistical Association\/}~{\em
  109\/}(508), 1601--1614.

\bibitem[\protect\citeauthoryear{Angrist and Pischke}{Angrist and
  Pischke}{2008}]{angristpischke}
Angrist, J. and S.~Pischke (2008).
\newblock {\em Mostly Harmless Econometrics: An Empiricists' Companion}.
\newblock Princeton University Press.

\bibitem[\protect\citeauthoryear{Angrist}{Angrist}{1998}]{angrist1998}
Angrist, J.~D. (1998).
\newblock Estimating the labor market impact of voluntary military service
  using social security data on military applicants.
\newblock {\em Econometrica\/}~{\em 66\/}(2), 249--288.

\bibitem[\protect\citeauthoryear{Aronow, Green, Lee, et~al.}{Aronow
  et~al.}{2014}]{aronow2014sharp}
Aronow, P.~M., D.~P. Green, D.~K. Lee, et~al. (2014).
\newblock Sharp bounds on the variance in randomized experiments.
\newblock {\em The Annals of Statistics\/}~{\em 42\/}(3), 850--871.

\bibitem[\protect\citeauthoryear{Aronow and Samii}{Aronow and
  Samii}{2016}]{aronow2016}
Aronow, P.~M. and C.~Samii (2016).
\newblock Does regression produce representative estimates of causal effects?
\newblock {\em American Journal of Political Science\/}~{\em 60\/}(1),
  250--267.

\bibitem[\protect\citeauthoryear{Davidson}{Davidson}{1994}]{davidson1994stochastic}
Davidson, J. (1994).
\newblock {\em Stochastic Limit Theory: An Introduction for Econometricians}.
\newblock Advanced Texts in Econometrics. Oxford University Press.

\bibitem[\protect\citeauthoryear{Deaton}{Deaton}{2010}]{deaton2010}
Deaton, A. (2010).
\newblock Instruments, randomization, and learning about development.
\newblock {\em Journal of Economic Literature\/}~{\em 48\/}(2), 424--455.

\bibitem[\protect\citeauthoryear{DellaVigna, Lindner, Reizer, and
  Schmieder}{DellaVigna et~al.}{2017}]{dellavigna2017reference}
DellaVigna, S., A.~Lindner, B.~Reizer, and J.~F. Schmieder (2017).
\newblock Reference-dependent job search: Evidence from hungary.
\newblock {\em The Quarterly Journal of Economics\/}~{\em 132\/}(4),
  1969--2018.

\bibitem[\protect\citeauthoryear{Donohue, Aneja, and Weber}{Donohue
  et~al.}{2017}]{donohue2017right}
Donohue, J.~J., A.~Aneja, and K.~Weber (2017).
\newblock Right-to-carry laws and violent crime: a comprehensive assessment
  using panel data, the lasso, and a state-level synthetic controls analysis.

\bibitem[\protect\citeauthoryear{Eicker}{Eicker}{1967}]{eicker}
Eicker, F. (1967).
\newblock Limit theorems for regressions with unequal and dependent errors.
\newblock In {\em Proceedings of the fifth Berkeley symposium on mathematical
  statistics and probability}, Volume~1, pp.\  59--82.

\bibitem[\protect\citeauthoryear{Einav, Finkelstein, and Schrimpf}{Einav
  et~al.}{2015}]{einav2015response}
Einav, L., A.~Finkelstein, and P.~Schrimpf (2015).
\newblock The response of drug expenditure to nonlinear contract design:
  evidence from medicare part d.
\newblock {\em The quarterly journal of economics\/}~{\em 130\/}(2), 841--899.

\bibitem[\protect\citeauthoryear{Farber}{Farber}{2015}]{farber2015you}
Farber, H.~S. (2015).
\newblock Why you can�t find a taxi in the rain and other labor supply
  lessons from cab drivers.
\newblock {\em The Quarterly Journal of Economics\/}~{\em 130\/}(4),
  1975--2026.

\bibitem[\protect\citeauthoryear{Fogarty}{Fogarty}{2016}]{fogarty2016regression}
Fogarty, C.~B. (2016).
\newblock Regression assisted inference for the average treatment effect in
  paired experiments.
\newblock {\em arXiv preprint arXiv:1612.05179\/}.

\bibitem[\protect\citeauthoryear{Freedman}{Freedman}{2008}]{Freedman2008}
Freedman, D. (2008).
\newblock On regression adjustmens to experimental data.
\newblock {\em Advances in Applied Mathematics\/}~{\em 40\/}(2), 180--193.

\bibitem[\protect\citeauthoryear{Goldberger}{Goldberger}{1991}]{goldberger1991course}
Goldberger, A.~S. (1991).
\newblock {\em A course in econometrics}.
\newblock Harvard University Press.

\bibitem[\protect\citeauthoryear{Hanna, Mullainathan, and Schwartzstein}{Hanna
  et~al.}{2014}]{hanna2014learning}
Hanna, R., S.~Mullainathan, and J.~Schwartzstein (2014).
\newblock Learning through noticing: Theory and evidence from a field
  experiment.
\newblock {\em The Quarterly Journal of Economics\/}~{\em 129\/}(3),
  1311--1353.

\bibitem[\protect\citeauthoryear{Holland}{Holland}{1986}]{holland1986statistics}
Holland, P.~W. (1986).
\newblock Statistics and causal inference.
\newblock {\em Journal of the American statistical Association\/}~{\em
  81\/}(396), 945--960.

\bibitem[\protect\citeauthoryear{Huber}{Huber}{1967}]{huber}
Huber, P.~J. (1967).
\newblock The behavior of maximum likelihood estimates under nonstandard
  conditions.
\newblock In {\em Proceedings of the fifth Berkeley symposium on mathematical
  statistics and probability}, Volume~1, pp.\  221--233.

\bibitem[\protect\citeauthoryear{Imbens and Rubin}{Imbens and
  Rubin}{2015}]{imbens2015causal}
Imbens, G.~W. and D.~B. Rubin (2015).
\newblock {\em Causal Inference in Statistics, Social, and Biomedical
  Sciences}.
\newblock Cambridge University Press.

\bibitem[\protect\citeauthoryear{Keels, Duncan, DeLuca, Mendenhall, and
  Rosenbaum}{Keels et~al.}{2005}]{keels2005fifteen}
Keels, M., G.~J. Duncan, S.~DeLuca, R.~Mendenhall, and J.~Rosenbaum (2005).
\newblock Fifteen years later: Can residential mobility programs provide a
  long-term escape from neighborhood segregation, crime, and poverty.
\newblock {\em Demography\/}~{\em 42\/}(1), 51--73.

\bibitem[\protect\citeauthoryear{Lin}{Lin}{2013}]{lin}
Lin, W. (2013).
\newblock Agnostic notes on regression adjustments for experimental data:
  Reexamining freedman's critique.
\newblock {\em The Annals of Applied Statistics\/}~{\em 7\/}(1), 295--318.

\bibitem[\protect\citeauthoryear{Lu}{Lu}{2016}]{lu2016randomization}
Lu, J. (2016).
\newblock On randomization-based and regression-based inferences for 2k
  factorial designs.
\newblock {\em Statistics \& Probability Letters\/}~{\em 112}, 72--78.

\bibitem[\protect\citeauthoryear{MacKinnon and White}{MacKinnon and
  White}{1985}]{mackinnon1985some}
MacKinnon, J.~G. and H.~White (1985).
\newblock Some heteroskedasticity-consistent covariance matrix estimators with
  improved finite sample properties.
\newblock {\em Journal of econometrics\/}~{\em 29\/}(3), 305--325.

\bibitem[\protect\citeauthoryear{Manski}{Manski}{2013}]{manski2013public}
Manski, C.~F. (2013).
\newblock {\em Public policy in an uncertain world: analysis and decisions}.
\newblock Harvard University Press.

\bibitem[\protect\citeauthoryear{Manski and Pepper}{Manski and
  Pepper}{2018}]{manski2018right}
Manski, C.~F. and J.~V. Pepper (2018).
\newblock How do right-to-carry laws affect crime rates? coping with ambiguity
  using bounded-variation assumptions.
\newblock {\em Review of Economics and Statistics\/}~{\em 100\/}(2), 232--244.

\bibitem[\protect\citeauthoryear{Munnell, Tootell, Browne, and
  McEneaney}{Munnell et~al.}{1996}]{munnell1996mortgage}
Munnell, A.~H., G.~M. Tootell, L.~E. Browne, and J.~McEneaney (1996).
\newblock Mortgage lending in boston: Interpreting hmda data.
\newblock {\em The American Economic Review\/}, 25--53.

\bibitem[\protect\citeauthoryear{Muralidharan and Niehaus}{Muralidharan and
  Niehaus}{2017}]{muralidharan2017experimentation}
Muralidharan, K. and P.~Niehaus (2017).
\newblock Experimentation at scale.
\newblock {\em Journal of Economic Perspectives\/}~{\em 31\/}(4), 103--24.

\bibitem[\protect\citeauthoryear{Neyman}{Neyman}{1990}]{neyman1923}
Neyman, J. (1923/1990).
\newblock On the application of probability theory to agricultural experiments.
  essay on principles. section 9.
\newblock {\em Statistical Science\/}~{\em 5\/}(4), 465--472.

\bibitem[\protect\citeauthoryear{Rosenbaum}{Rosenbaum}{2002}]{rosenbaum_book}
Rosenbaum, P.~R. (2002).
\newblock {\em Observational Studies}.
\newblock Springer.

\bibitem[\protect\citeauthoryear{Rosenbaum and Rubin}{Rosenbaum and
  Rubin}{1983}]{rosenbaum1983central}
Rosenbaum, P.~R. and D.~B. Rubin (1983).
\newblock The central role of the propensity score in observational studies for
  causal effects.
\newblock {\em Biometrika\/}~{\em 70\/}(1), 41--55.

\bibitem[\protect\citeauthoryear{Rubin}{Rubin}{1974}]{rubin1974estimating}
Rubin, D.~B. (1974).
\newblock Estimating causal effects of treatments in randomized and
  nonrandomized studies.
\newblock {\em Journal of educational Psychology\/}~{\em 66\/}(5), 688.

\bibitem[\protect\citeauthoryear{Samii and Aronow}{Samii and
  Aronow}{2012}]{samii2012equivalencies}
Samii, C. and P.~M. Aronow (2012).
\newblock On equivalencies between design-based and regression-based variance
  estimators for randomized experiments.
\newblock {\em Statistics \& Probability Letters\/}~{\em 82\/}(2), 365--370.

\bibitem[\protect\citeauthoryear{Shadish, Cook, and Campbell}{Shadish
  et~al.}{2002}]{shadishcookcampbell}
Shadish, W.~R., T.~D. Cook, and D.~T. Campbell (2002).
\newblock {\em Experimental and quasi-experimental designs for generalized
  causal inference.}
\newblock Houghton, Mifflin and Company.

\bibitem[\protect\citeauthoryear{S\l{}oczy\'{n}ski}{S\l{}oczy\'{n}ski}{2017}]{tymon2017}
S\l{}oczy\'{n}ski, T. (2017).
\newblock A general weighted average representation of the ordinary and
  two-stage least squares estimands.

\bibitem[\protect\citeauthoryear{White}{White}{1980a}]{white1980robust}
White, H. (1980a).
\newblock A heteroskedasticity-consistent covariance matrix estimator and a
  direct test for heteroskedasticity.
\newblock {\em Econometrica\/}~{\em 48\/}(1), 817--838.

\bibitem[\protect\citeauthoryear{White}{White}{1980b}]{white1980using}
White, H. (1980b).
\newblock Using least squares to approximate unknown regression functions.
\newblock {\em International Economic Review\/}~{\em 21\/}(1), 149--170.

\bibitem[\protect\citeauthoryear{White}{White}{1982}]{white1982maximum}
White, H. (1982).
\newblock Maximum likelihood estimation of misspecified models.
\newblock {\em Econometrica\/}~{\em 50\/}(1), 1--25.

\bibitem[\protect\citeauthoryear{Wooldridge}{Wooldridge}{2001}]{wooldridge2001}
Wooldridge, J.~M. (2001).
\newblock Asymptotic properties of weighted {M}-estimators for standard
  stratified samples.
\newblock {\em Econometric Theory\/}~{\em 17\/}(2), 451--470.

\end{thebibliography}

\newpage

\begin{table}[ht]
\caption{\textsc{: Simulation Results  with coverage
for nominal 95\% confidence intervals}}
\label{tabel1}\vskip1cm
\par
\begin{center}
\begin{tabular}{lccccccccc}
$E[N]=\rho n$& 1000& 1000& 100& 1000& 1000& 1000& 1000\\
$\rho$ & 0.01& 0.01& 0.01& 1& 0.01& 0.01& 0.01\\
$k$& 1& 10& 1& 1& 1& 1& 1\\
$\psi'\psi$&4 &4&4&4&0&4&0\\
$\sigma^2_\theta$ & 1 & 1 & 1 & 1 & 1 & 0 & 0\\[2ex]

sd($\widehat\theta_n- \theta_n^\mathrm{descr}$)                       			& 0.125 & 0.126 & 0.399 & 0.000 & 0.063 & 0.113 & 0.031 \\
sd($\widehat\theta_n- \theta_n^{\mathrm{causal},\mathrm{sample}}$)    	& 0.105 & 0.104 & 0.331 & 0.100 & 0.055 & 0.095 & 0.032 \\
sd($\widehat\theta_n- \theta_n^{\mathrm{causal}}$)                    			& 0.125 & 0.126 & 0.400 & 0.100 & 0.063 & 0.114 & 0.032 \\[2ex]

average $(\smash[b]{\widehat V_n^{\mathrm{ehw}}/N})^{1/2}$            	& 0.125 & 0.124 & 0.370 & 0.121 & 0.063 & 0.113 & 0.032 \\
coverage $\theta_n^{\mathrm{descr}}$                                  				& 0.949 & 0.947 & 0.923 & 1.000 & 0.948 & 0.947 & 0.950 \\
coverage $\theta_n^{\mathrm{causal},\mathrm{sample}}$                 		& 0.980 & 0.981 & 0.969 & 0.982 & 0.974 & 0.981 & 0.950 \\
coverage $\theta_n^{\mathrm{causal}}$                                 				& 0.948 & 0.947 & 0.922 & 0.982 & 0.947 & 0.947 & 0.950 \\[2ex]

average $(\smash[b]{\widehat V_n^{\mathrm{boot}}/N})^{1/2}$           	& 0.126 & 0.127 & 0.426 & 0.122 & 0.064 & 0.115 & 0.032 \\
coverage $\theta_n^{\mathrm{descr}}$                                  				& 0.950 & 0.950 & 0.955 & 1.000 & 0.950 & 0.949 & 0.953 \\
coverage $\theta_n^{\mathrm{causal},\mathrm{sample}}$                 		& 0.981 & 0.982 & 0.986 & 0.982 & 0.975 & 0.981 & 0.951 \\
coverage $\theta_n^{\mathrm{causal}}$                                 				& 0.950 & 0.949 & 0.955 & 0.982 & 0.949 & 0.949 & 0.951  \\[2ex]

average $(\smash[b]{\widehat V_n^{\mathrm{desc}}/N})^{1/2}$           	& 0.124 & 0.124 & 0.368 & 0.000 & 0.063 & 0.113 & 0.031 \\
coverage $\theta_n^{\mathrm{descr}}$                                  				& 0.948 & 0.946 & 0.921 & 1.000 & 0.947 & 0.946 & 0.949 \\
coverage $\theta_n^{\mathrm{causal},\mathrm{sample}}$                 		& 0.980 & 0.980 & 0.968 & 0.000 & 0.973 & 0.981 & 0.948 \\
coverage $\theta_n^{\mathrm{causal}}$                                 				& 0.947 & 0.946 & 0.921 & 0.000 & 0.946 & 0.946 & 0.948 \\[2ex]

average $(\smash[b]{\widehat V_n^{\causal,\sample}/N})^{1/2}$         	& 0.108 & 0.107 & 0.317 & 0.104 & 0.063 & 0.094 & 0.032 \\
coverage $\theta_n^{\mathrm{descr}}$                                  				& 0.908 & 0.905 & 0.872 & 1.000 & 0.948 & 0.894 & 0.950 \\
coverage $\theta_n^{\mathrm{causal},\mathrm{sample}}$                 		& 0.956 & 0.957 & 0.937 & 0.957 & 0.974 & 0.948 & 0.949 \\
coverage $\theta_n^{\mathrm{causal}}$                                 				& 0.907 & 0.904 & 0.870 & 0.957 & 0.947 & 0.892 & 0.949 \\[2ex]

average $(\smash[b]{\widehat V_n^{\causal}/N})^{1/2}$                 		& 0.125 & 0.124 & 0.369 & 0.104 & 0.063 & 0.113 & 0.032 \\
coverage $\theta_n^{\mathrm{descr}}$                                  				& 0.949 & 0.947 & 0.922 & 1.000 & 0.948 & 0.947 & 0.950 \\
coverage $\theta_n^{\mathrm{causal},\mathrm{sample}}$                 		& 0.980 & 0.981 & 0.969 & 0.957 & 0.974 & 0.981 & 0.950 \\
coverage $\theta_n^{\mathrm{causal}}$                                 				& 0.948 & 0.947 & 0.922 & 0.957 & 0.947 & 0.946 & 0.950
\end{tabular}%
\end{center}
\end{table}

\end{document}